\documentclass[12pt,a4paper]{article}

\usepackage{amsmath,amssymb,amsbsy,amsfonts,amsthm,latexsym,
mathabx,amsopn,amstext,amsxtra,euscript,amscd,stmaryrd,mathrsfs,cite,array,mathtools,enumerate,url}

\usepackage[colorlinks,linkcolor=blue,anchorcolor=blue,citecolor=blue,backref=page]{hyperref}

\usepackage{geometry}
\usepackage[dvipsnames]{xcolor}

\usepackage{soul}

\usepackage{pdfpages}
\usepackage{fancyhdr} 
\usepackage{bookmark}

\newtheorem{theorem}{Theorem}
\newtheorem{lemma}[theorem]{Lemma}

\newtheorem{cor}[theorem]{Corollary}

\theoremstyle{remark}

\newtheorem*{example}{Example}

\theoremstyle{definition}
\newtheorem{definition}{Definition}

\newtheorem{problem}{Problem}

\numberwithin{theorem}{section} 
\numberwithin{definition}{section} 

\newcommand{\F}{\mathbb{F}}
\newcommand{\N}{\mathbb{N}}

\newcommand{\R}{\mathbb{R}}

\def\cA{{\mathcal A}}

\def\cL{{\mathcal L}}

\usepackage{tikz}

\usepackage{todonotes}

\numberwithin{equation}{section}

\providecommand{\keywords}[1]
{
  \small \noindent	
  {\textit{Keywords ---}} #1
}

\providecommand{\MSC}[1]
{
  \small	\noindent	
  {\textit{2020 Mathematics Subject Classification
---}} #1
}

\title{Pseudorandom sequences derived from automatic sequences}
\date{}

\author{L\'aszl\'o  M\'erai and Arne Winterhof\\
\small
Johann Radon Institute for Computational and Applied Mathematics \\
\small Austrian Academy of Sciences,
Altenbergerstr.\ 69, 4040 Linz, Austria\\
\small E-mail: \texttt{\{laszlo.merai,arne.winterhof\}@oeaw.ac.at}}

\begin{document}

\maketitle


\begin{abstract}
    Many automatic sequences, such as the Thue-Morse sequence or the Rudin-Shapiro sequence,
    have some desirable features of pseudorandomness such as a large linear complexity and a small well-distribution measure.
    However, they also have some disastrous properties in view of certain applications. For example, the majority of possible binary patterns never appears in automatic sequences and their correlation measure of order $2$ is extremely large.
    
    Certain subsequences, such as  automatic sequences along squares, may keep the good properties of the original sequence but avoid the bad ones.  
    
     In this survey we investigate properties of pseudorandomness and non-randomness of automatic sequences and their subsequences and present results on their behaviour under several measures of pseudorandomness including linear complexity, correlation measure of order $k$, expansion complexity and normality.
     We also mention some analogs for finite fields.
\end{abstract}

\keywords{automatic sequences, pseudorandomness, linear complexity, maximum order complexity, well-distribution measure, correlation measure, expansion complexity, normality, finite fields}
\smallskip

\MSC{11A63, 11B85, 11K16, 11K31, 11K36, 11K45, 11T71, 68R15, 94A55, 94A60}

\section{Introduction}

Pseudorandom sequences are sequences generated by deterministic algorithms which shall simulate randomness. 
In contrast to truly random sequences they are not random at all but guarantee certain desirable features and are reproducible.

Automatic sequences, see Section~\ref{sec:automatic_sequences} below for the definition, have some of these desirable features but also some undesirable ones.

\begin{figure}
\begin{center}
\includegraphics[scale=0.15]{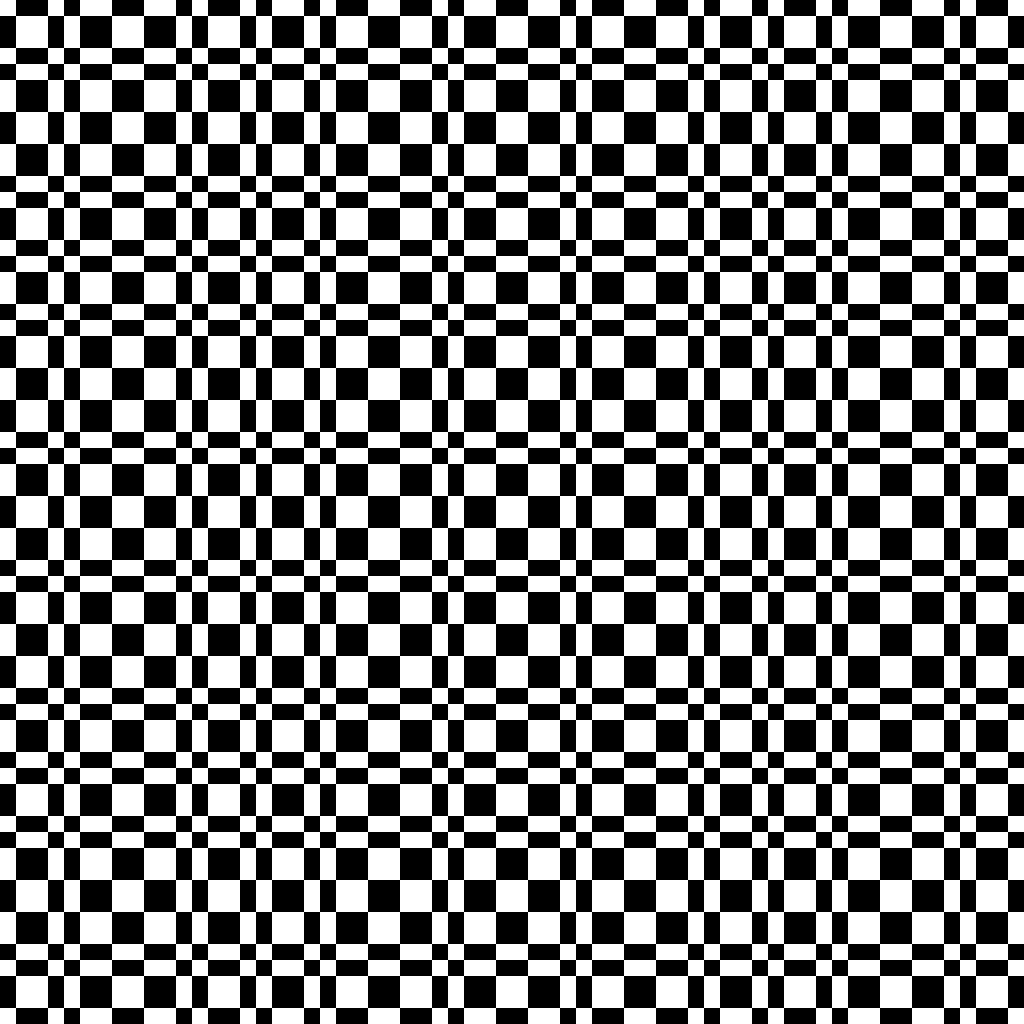} \qquad
\includegraphics[scale=0.15]{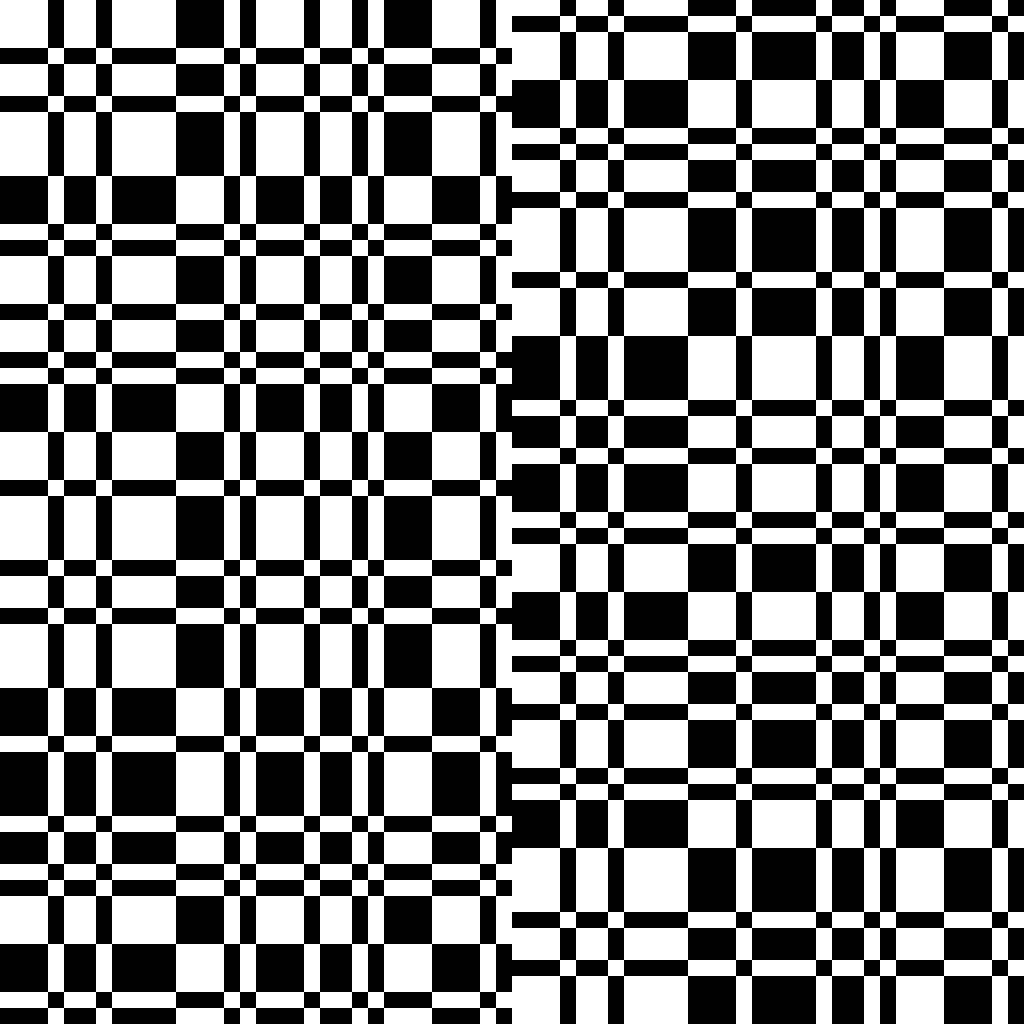} 
\end{center}
\caption{The first $4096$ elements of the Thue-Morse (left) and Rudin-Shapiro (right) sequence  split into $64$ rows of each $64$ sequence elements. Zeros are represented by white, ones are represented by black.}
\label{tmfig}
\end{figure}

For example, the Thue-Morse sequence $(t_n)$, defined by \eqref{tmdef} below, 
\begin{itemize}
\item has large $N$th linear complexity, see Section~\ref{sec:linear_complexity}, 
\item has large $N$th maximum-order complexity, see Section~\ref{sec:max-order_complexity},
\item is balanced and has a small well-distribution measure, see Section~\ref{sec:correlation}.
\end{itemize}
However, the Thue-Morse sequence 
\begin{itemize}
\item has a very large correlation measure of order $2$,
see Section~\ref{sec:correlation}, 
\item a very small expansion complexity, see Section~\ref{sec:expansion_complexity}, 
\item and there are short patterns such as $000$ and $111$ which do not appear in the sequence and its subword complexity is only linear, see Section~\ref{sec:normality}.
\end{itemize}
Hence, despite 
some nice features this sequence is not looking random at all, see Figure~\ref{tmfig}. The same is true for the Rudin-Shapiro sequence $(r_n)$ defined by~\eqref{rsdef} below and many other related sequences.

\begin{figure}[ht] 
\begin{center}
\includegraphics[scale=0.15]{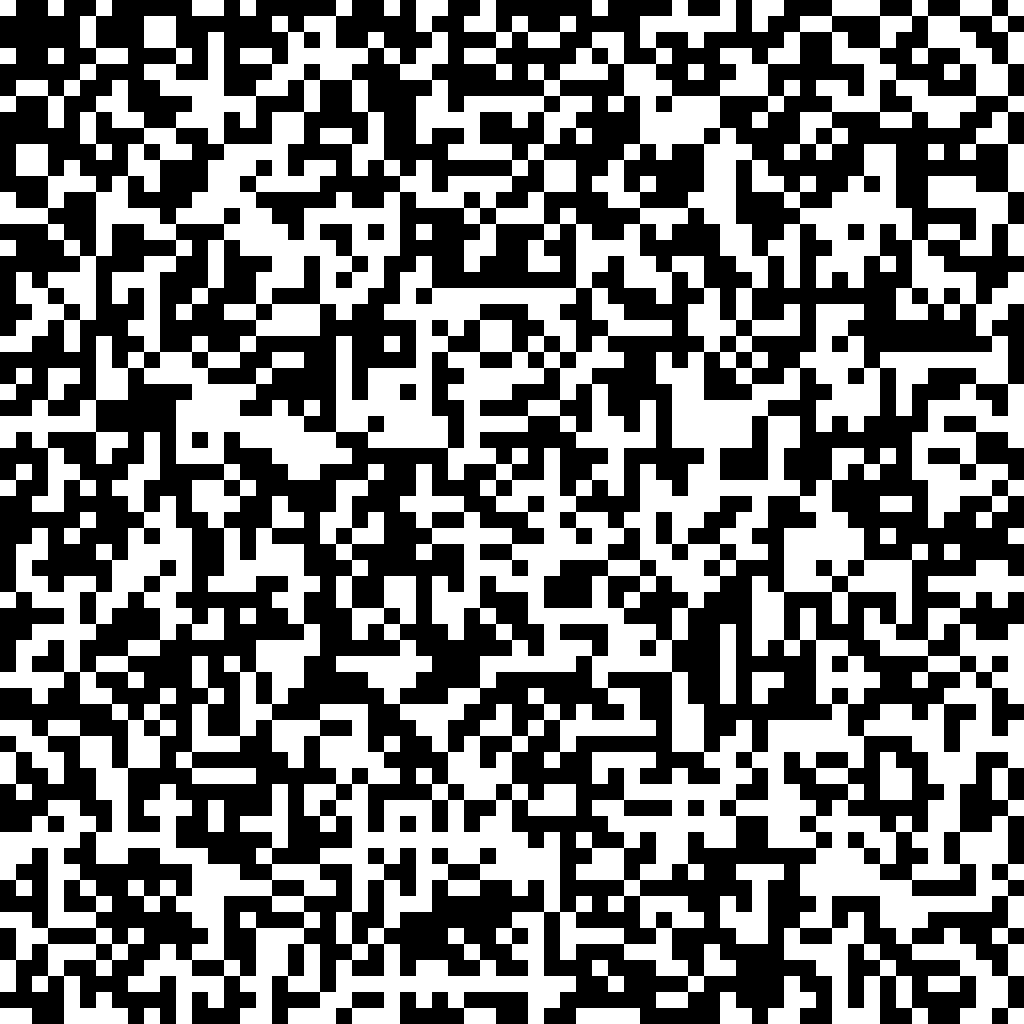} \qquad
\includegraphics[scale=0.15]{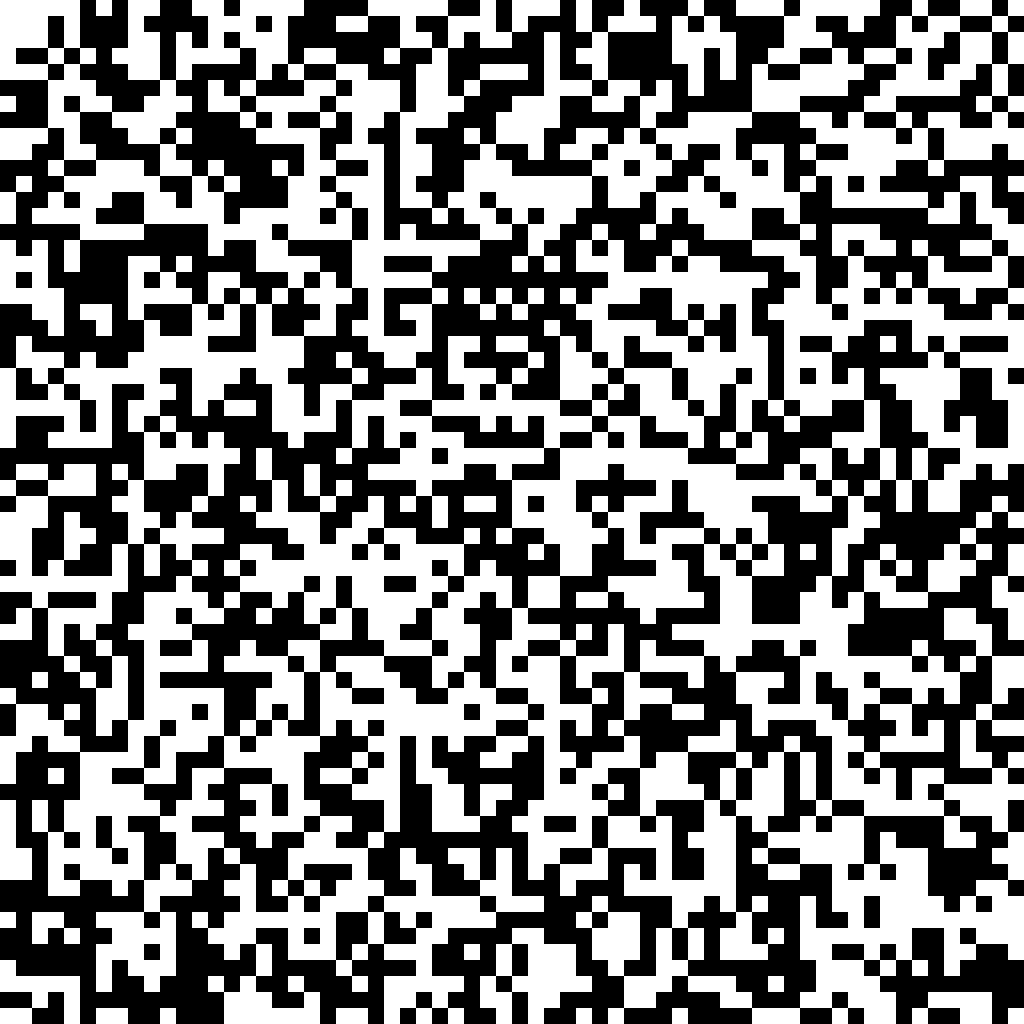} 
\end{center}
\caption{The first $4096$ elements of the Thue-Morse (left) and Rudin-Shapiro (right) sequence along squares
split into $64$ rows of each $64$ sequence elements. Zeros are represented by white, ones are represented by black.}
\label{squarefig}
\end{figure}

Taking suitable subsequences may destroy the non-random structure of the original sequence but may keep the desirable features of pseudorandomness.
Promising candidates for such subsequences are
\begin{itemize}
\item along squares, cubes, bi-squares, ... or any polynomial values
    for any polynomial $f$ of degree at least $2$ with $f(\N_0)\subset \N_0$,
\item along primes,
\item along the Piateski-Shapiro sequence $\lfloor n^c\rfloor$, $1<c<2$,
\item and along geometric sequences such as $3^n$.
\end{itemize}
For example, the Thue-Morse sequence and the Rudin-Shapiro sequence along squares still 
\begin{itemize}
    \item have a large maximum-order complexity and thus a large linear complexity, see Section~\ref{sec:max-order_complexity},
    \item and are asymptotically balanced, see Section~\ref{sec:normality}. 
\end{itemize}
Moreover, in contrast to the original sequence they 
\begin{itemize}
    \item have unbounded expansion complexity, see Section~\ref{sec:expansion_complexity},
    \item and are normal, that is, asymptotically each pattern appears with the right frequency in the sequence, see Section~\ref{sec:normality}.
\end{itemize}
Roughly speaking, they look much more random than the original sequences, see Figure~\ref{squarefig}.

Still some questions about these sequences remain open such as upper bounds on the correlation measure of order $k$ and on the expansion complexity. We will state explicitly some selected open problems to motivate future research.

We also look for further directions in Section~\ref{sec:finite_fields}.  In particular, we discuss analogs of the Thue-Morse and Rudin-Shapiro sequence and their subsequences in the setting of finite fields. 

For general background on automatic sequences and finite automata we refer to the monograph of Allouche and Shallit \cite{alsh2} and also to \cite{alsh,alshya,ev,fo}.
For surveys on  pseudorandom sequences see \cite{gy,merisa,niwi,sh,towi}.

\section{Finite automata and automatic sequences}\label{sec:automatic_sequences}

Roughly speaking, a sequence is {\em automatic} if it is  generated by a finite automaton, see Definition~\ref{def:sequence} below.

\begin{definition}
Let $k\geq 2$ be an integer.
A \emph{finite $k$-automaton} $\cA$ is a $6$-tuple 
$$
\cA=(Q,\Sigma, \delta, q_0, \varphi,\Delta),
$$ 
where
\begin{itemize}
    \item $Q$ is a finite set of states, 
    \item $\Sigma=\{0,1,\ldots,k-1\}$ is the input alphabet,  
    \item $\delta: Q \times \Sigma \rightarrow Q$ is the transition function, 
    \item $q_0\in Q$ is the initial state, 
    \item $\Delta$ is the output alphabet
    \item and $\varphi:Q \rightarrow \Delta$ is the output function.
\end{itemize}
\end{definition}

\begin{figure}[ht]
\begin{center}
\begin{tikzpicture}[auto,thick]
 \node (E) at (-2,0) [circle] {};
 \node (A) at (0,0) [circle, draw] {$A/0$};
 \node (B) at (3,0) [circle, draw] {$B/0$};
 \node (C) at (6,0) [circle, draw] {$C/1$};
 \node (D) at (9,0) [circle, draw] {$D/1$};

 \draw [->,bend left] (A) to  node {1} (B);
 \draw [->,bend left] (B) to  node {0} (A);
\draw [->,bend left] (B) to  node {1} (C);
 \draw [->,bend left] (C) to  node {1} (B);
\draw [->,bend left] (C) to  node {0} (D);
 \draw [->,bend left] (D) to  node {1} (C);

 \path (D) edge [loop above] node {0} (B);
 \path (A) edge [loop above] node {0} (A);

 \draw [->] (E) to  node 
 [midway,above,align=center ] {\texttt{start}}
 (A);
\end{tikzpicture}
\end{center}
\caption{Rudin-Shapiro automaton} \label{fig:RS}
\end{figure}
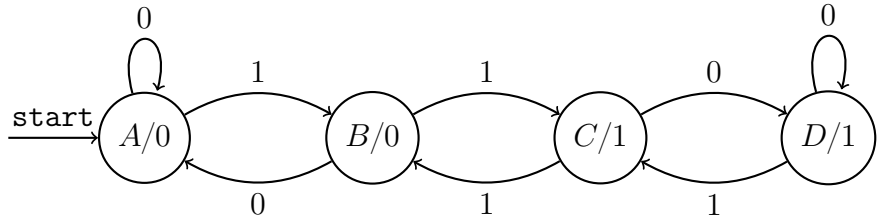

For example, the \emph{Thue-Morse automaton}, see Figure~\ref{fig:TM}, is a $2$-automaton with $2$ states and the \emph{Rudin-Shapiro automaton}, see Figure~\ref{fig:RS}, is a $2$-automaton with $4$ states, both with inputs and outputs in $\Sigma=\Delta=\{0,1\}$. 

\begin{figure}[ht]
\begin{center}
\begin{tikzpicture}[auto,thick]
 \node (E) at (-2,0) [circle] {};
 \node (A) at (0,0) [circle, draw] {$A/0$};
 \node (B) at (3,0) [circle, draw] {$B/1$};
 \draw [->,bend left] (A) to  node {1} (B);
 \draw [->,bend left] (B) to  node {1} (A);
 \path (B) edge [loop above] node {0} (B);
 \path (A) edge [loop above] node {0} (A);
 \draw [->] (E) to  node 
 [midway,above,align=center ] {\texttt{start}}
 (A);
\end{tikzpicture}
\end{center}
\caption{Thue-Morse automaton} \label{fig:TM}
\end{figure}
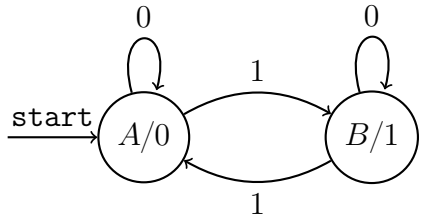

\begin{definition}\label{def:sequence}
 Let $\Delta$ be a finite set. A sequence $(s_n)$ over $\Delta$ is called a \emph{$k$-automatic sequence} if there is a  $k$-automaton $\cA$ such that on input of the digits $n_0,n_1,\ldots$ of the $k$-ary expansion of $n\geq 0$, 
\begin{equation}\label{eq:k-ary}
 n=\sum_{i\geq 0} n_i k^i, \quad n_i\in\{0,1,\dots, k-1\},
\end{equation}
 $\cA$ outputs the sequence element $s_n\in \Delta$. 
Reading of the digits of $n$ starting with the most significant digit is called {\em direct}  whereas reading starting with the least significant digit $n_0$ is called {\em reverse}.
If not stated otherwise, we use reverse reading.
Finally, a sequence is called {\em automatic} if it is $k$-automatic for some $k$.
\end{definition}

\begin{example}[Thue-Morse sequence]
The \emph{Thue-Morse sequence}~$(t_n)$ is a $2$-automatic sequence generated by the Thue-Morse automaton, Figure~\ref{fig:TM}. This sequence is the {\em sequence of the sum of digits modulo $2$}. The sequence begins with
$$
011010011001 \dots,
$$
see also Figure~\ref{tmfig} for a picture of the first $4096$ sequence elements.
It follows from the defining automaton, see Figure~\ref{fig:TM}, that $(t_n)$
satisfies the following recurrence relation
\begin{equation}\label{tmdef}
t_n=
\left\{ 
\begin{array}{cl} 
t_{n/2} & \mbox{if $n$ is even},\\ 
t_{(n-1)/2}+1 \bmod 2 & \mbox{if $n$ is odd},
\end{array}\right.
\quad n=1,2,\ldots
\end{equation}
with initial value $t_0=0$.
\end{example}

\begin{example}[Rudin-Shapiro sequence]
The \emph{Rudin-Shapiro sequence} $(r_n)$ is a $2$-automatic sequence generated by the Rudin-Shapiro automaton, see Figure~\ref{fig:RS}. The sequence begins with
$$
000100100001   \dots,
$$
see also Figure~\ref{tmfig} for a picture of the first $4096$ sequence elements.
It follows from the defining automaton, see Figure~\ref{fig:RS}, that $(r_n)$
satisfies the following recurrence relation
\begin{equation}\label{rsdef}
r_n=
\left\{ 
\begin{array}{cl} 
r_{\lfloor n/2\rfloor}+1 \bmod 2 & \mbox{if $n\equiv 3 \bmod 4$},\\ 
r_{\lfloor n/2\rfloor} & \mbox{otherwise},
\end{array}\right.
\quad n=1,2,\ldots
\end{equation}
with initial value $r_0=0$.
\end{example}
The sequence $((-1)^{r_n})$ over $\{-1,+1\}$ is also called Rudin-Shapiro sequence in the literature. Here we study only the sequence $(r_n)$ over $\{0,1\}$.

\begin{example}[Pattern sequences]
For a pattern 
$P\in \Delta^\ell\setminus\{(0,\dots,0)\}$ of length $\ell$
over $\Delta=\{0,1,\dots, k-1\}$
define the sequence $(p_n)$ by
\begin{equation*}
 p_n= e_P(n) \bmod k, \quad 0\leq p_n<k,  \quad n=0,1,\dots,
\end{equation*}
where $e_P(n)$ is the number of occurrences of $P$ in the $k$-ary expansion of $n$. 
The sequence~$(p_n)$ over~$\Delta$ satisfies the following recurrence relation
\begin{equation}\label{eq:recurrence}
 p_n=
 \left\{
 \begin{array}{cl}
 p_{\lfloor n/k\rfloor}+1 \bmod k  & \text{if } n\equiv a \bmod k^\ell,\\
 p_{\lfloor n/k\rfloor}  & \text{otherwise,}
 \end{array}
 \right.
 n=1,2,\dots
\end{equation}
with initial value $p_0=0$, where $a=a(P)$ is the integer $0< a <k^\ell$ such that its  $k$-ary expansion corresponds to the pattern $P$.

Classical examples for binary pattern sequences are the Thue-Morse sequence with 
$$k=2,\quad \ell=1, \quad P=1 \quad \mbox{and}\quad a=1,$$ and the Rudin-Shapiro sequence with 
$$k=2,\quad \ell=2,\quad  P=11\quad \mbox{and}\quad a=3.$$
In particular, if $n_0,n_1,\dots$ are the bits of the non-negative integer taken from \eqref{eq:k-ary} with $k=2$,
then
\begin{equation}\label{sumofdigitsdef}
t_n=\sum_{i=0}^\infty n_i \bmod 2 \quad \mbox{and}\quad r_n=\sum_{i=0}^\infty n_in_{i+1} \bmod 2.
\end{equation}
\end{example}

\begin{example}[Rudin-Shapiro-like sequence]
Lafrance, Rampersad and Yee \cite{laraye} introduced a {\em Rudin-Shapiro-like sequence}
$(\ell_n)$ which is based on the number of occurrences of the pattern $10$ as a scattered subsequence in the binary representation, \eqref{eq:k-ary} with $k=2$, of 
$n$. That is,~$\ell_n$ is
the  parity of the number of pairs $(i,j)$ with $i>j$ and $(n_i,n_j)=(1,0)$. See Figure~\ref{fig:RS-like} for its defining automaton.
\begin{figure}[ht]
\begin{center}
\begin{tikzpicture}[auto,thick]
 \node (E) at (-2,0) [circle] {};
 \node (A) at (0,0) [circle, draw] {$A/0$};
 \node (B) at (3,0) [circle, draw] {$B/0$};
 \node (C) at (6,0) [circle, draw] {$C/1$};
 \node (D) at (9,0) [circle, draw] {$D/1$};

 \draw [->,bend left] (A) to  node {1} (B);
 \draw [->,bend left] (B) to  node {1} (A);
\draw [->,bend left] (B) to  node {0} (C);
 \draw [->,bend left] (C) to  node {0} (B);
\draw [->,bend left] (C) to  node {1} (D);
 \draw [->,bend left] (D) to  node {1} (C);

 \path (D) edge [loop above] node {0} (B);
 \path (A) edge [loop above] node {0} (A);

 \draw [->] (E) to  node 
 [midway,above,align=center ] {\texttt{start}}
 (A);
\end{tikzpicture}

\end{center}
\caption{Rudin-Shapiro-like automaton with direct reading} \label{fig:RS-like}
\end{figure}
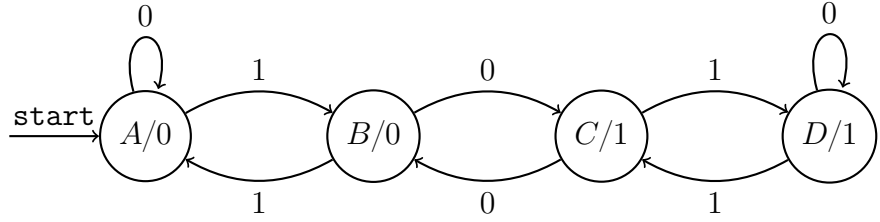
This sequence can also be defined by
\begin{equation}\label{eq:rslike}\ell_{2n+1}=\ell_n \quad \mbox{and}\quad \ell_{2n}=\ell_n+t_n\bmod 2, 
\end{equation}
see \cite[$(1)$ and $(2)$]{laraye},
with initial value $\ell_0=0$ and where $(t_n)$ is the Thue-Morse sequence.
\end{example}

\begin{figure}[ht]
\begin{center}
\begin{tikzpicture}[auto,thick]
 \node (E) at (-2,0) [circle] {};
 \node (A) at (0,0) [circle, draw] {$A/1$};
 \node (B) at (3,0) [circle, draw] {$B/1$};
 \node (C) at (6,0) [circle, draw] {$C/0$};

 \path (A) edge [loop above] node {1} (A);
 \draw [->] (E) to  node 
 [midway,above,align=center ] {\texttt{start}}
 (A);

 \draw [->,bend left] (A) to  node {0} (B);
 
 \draw [->,bend left] (B) to  node {0} (A);
\draw [->,left] (B) to  node [midway,above,align=center ] {1} (C);

 \path (C) edge [loop above] node {0,1} (C);

\end{tikzpicture}
\end{center}
\caption{Baum-Sweet automaton} \label{fig:BS}
\end{figure}
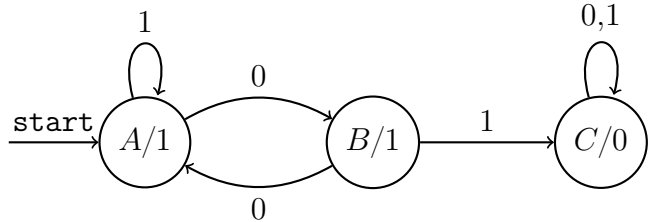

\begin{example}[Baum-Sweet sequence]
The \emph{Baum-Sweet sequence} $(b_n)$ is a $2$-automatic sequence defined by the rule $b_0=1$ and for $n\ge 1$
$$
 b_n=
 \left\{ 
 \begin{array}{cl}
  1& \text{if the binary representation of $n$ contains no block of} \\
   &  \text{consecutive $0$'s of odd length,}\\
  0& \text{otherwise.}
 \end{array}
 \right.
$$
Equivalently, we have for $n\geq 1$ of the form $n=4^\ell m$ with $4\nmid m$  that
\begin{equation}\label{bsdef}
 b_n=
 \left\{ 
 \begin{array}{cl}
  0& \text{if $m$ is even}, \\
  b_{(m-1)/2}& \text{if $m$ is odd.}
 \end{array}
 \right.
\end{equation}
The sequence 
$(b_n)$ is generated by the Baum-Sweet automaton in Figure~\ref{fig:BS}. 
\end{example}

\begin{example}[Characteristic sequence of sums of three squares] 
Consider the {\em characteristic sequence $(c_n)$ of the set of integers which are sums of three squares of an integer},
that is, 
$$ 
c_n=
 \left\{ 
 \begin{array}{cl}
  1& \text{if } n=a^2+b^2+c^2  \text{ for some non-negative integers $a,b,c$,} \\
  0& \text{otherwise.}
 \end{array}
 \right.
$$
By Legendre's three-square theorem, we have the equivalent definition
\begin{equation}\label{eq:cndef}
 c_n=
 \left\{ 
 \begin{array}{cl}
  1& \text{if $n$ is not of the form $n=4^\ell(8k+7)$,} \\
  0& \text{otherwise.}
 \end{array}
 \right.
\end{equation}
See Figure~\ref{fig:3squares} for the defining automaton. 

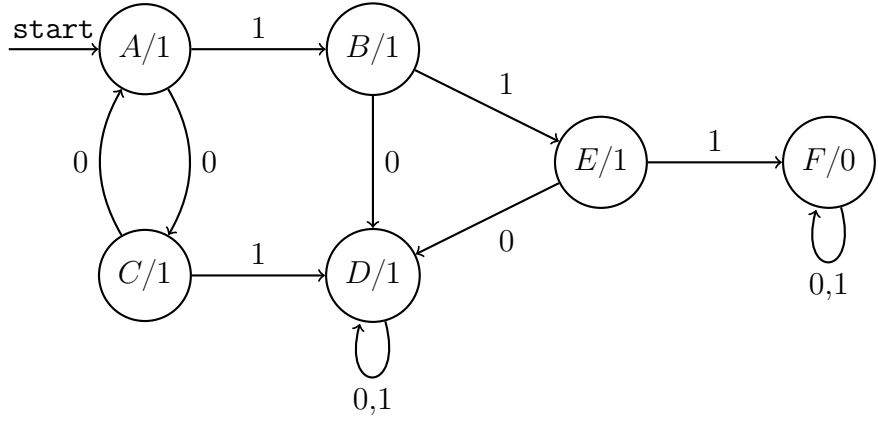
\begin{figure}[ht]
\begin{center}
\begin{tikzpicture}[auto,thick]
 \node (AA) at (-2,0) [circle] {};
 \node (A) at (0,0) [circle, draw] {$A/1$};
 \node (B) at (3,0) [circle, draw] {$B/1$};
 \node (C) at (0,-3) [circle, draw] {$C/1$};
 \node (D) at (3,-3) [circle, draw] {$D/1$};
 \node (E) at (6,-1.5) [circle, draw] {$E/1$};
 \node (F) at (9,-1.5) [circle, draw] {$F/0$};

 \draw [->] (A) to  node {1} (B);
 \draw [->,bend left] (A) to  node {0} (C);
 
 \draw [->] (B) to  node {1} (E);
 \draw [->] (B) to  node {0} (D);

 \draw [->] (C) to  node {1} (D);
 \draw [->,bend left] (C) to  node {0} (A);

 \draw [->] (E) to  node {1} (F);
 \draw [->] (E) to  node {0} (D);

 \path (D) edge [loop below] node {0,1} (D);
 \path (F) edge [loop below] node {0,1} (F);

 \draw [->] (AA) to  node 
 [midway,above,align=center ] {\texttt{start}}
 (A);
\end{tikzpicture}

\end{center}
\caption{Automaton of the characteristic sequence of sums of three squares with reverse reading} \label{fig:3squares}
\end{figure}
\end{example}

\begin{example}[Regular paper-folding sequence]
The {\em regular paper-folding sequence} $(v_n)$
with initial value $v_0\in \{0,1\}$ is defined as follows. If $n=2^km$ with an odd $m$, then
\begin{equation}\label{pfdef}
v_n=\left\{\begin{array}{ll}1,& m\equiv 1\bmod 4,\\ 0, &m\equiv 3\bmod 4,\end{array}\right.\quad n=1,2,\ldots
\end{equation}
Its defining automaton with four states is given in Figure~\ref{fig:RPF}.
\end{example}
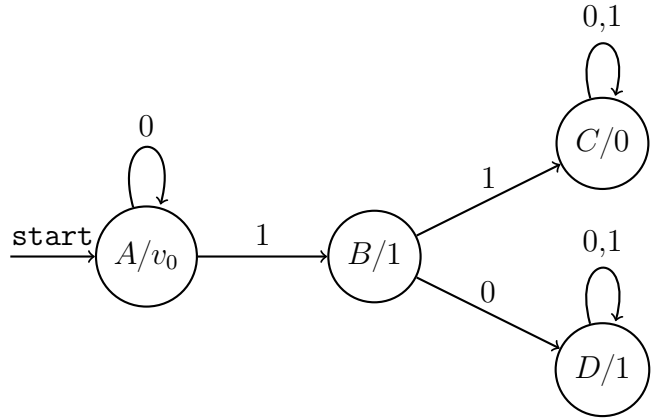
\begin{figure}[ht]
\begin{center}
\begin{tikzpicture}[auto,thick]
 \node (E) at (-2,0) [circle] {};
 \node (A) at (0,0) [circle, draw] {$A/v_0$};
 \node (B) at (3,0) [circle, draw] {$B/1$};
 \node (C) at (6,1.5) [circle, draw] {$C/0$};
 \node (D) at (6,-1.5) [circle, draw] {$D/1$};

 \path (A) edge [loop above] node {0} (A);
 \draw [->] (E) to  node 
 [midway,above,align=center ] {\texttt{start}}
 (A);

 \draw [->] (A) to  node {1} (B);
 
 \draw [->] (B) to  node [midway,above,align=center ]  {1} (C);
 \draw [->] (B) to  node [midway,above,align=center ]  {0} (D);
 
 \path (C) edge [loop above] node {0,1} (C);
 \path (D) edge [loop above] node {0,1} (D);

\end{tikzpicture}
\end{center}
\caption{Regular paper-folding automaton} \label{fig:RPF}
\end{figure}

\begin{example}[An automatic apwenian sequence]
Any binary sequence $(a_n)$ satisfying $a_0=1$
and 
$$
a_{2n+2}=a_{2n+1}+a_n \bmod 2,\quad n=0,1,\ldots
$$
is called {\em apwenian}, see for example \cite{alhani}. Apwenian sequences which are $2$-automatic are characterized in \cite{alhani}. For example, the sequence $(w_n)$ defined by 
\begin{equation}\label{apdef}
w_{2n}=1\quad \mbox{and}\quad w_{2n+1}=w_n+1\bmod 2,\quad n=0,1,\ldots
\end{equation}
is apwenian and defined by the automaton in Figure~\ref{fig:apw}.
\end{example}
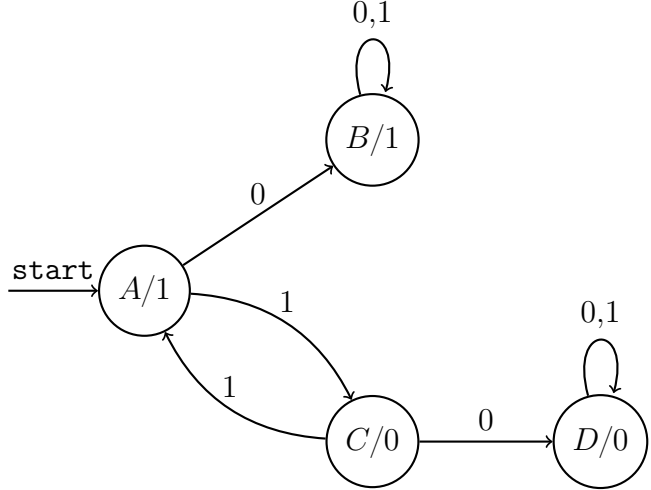
\begin{figure}[ht]
\begin{center}
\begin{tikzpicture}[auto,thick]
 \node (E) at (-2,0) [circle] {};
 \node (A) at (0,0) [circle, draw] {$A/1$};
 \node (B) at (3,2) [circle, draw] {$B/1$};
 \node (C) at (3,-2) [circle, draw] {$C/0$};
 \node (D) at (6,-2) [circle, draw] {$D/0$};

 \draw [->] (E) to  node 
 [midway,above,align=center ] {\texttt{start}}
 (A);

 \draw [->] (A) to  node [midway,above,align=center ]  {0} (B);
 
 \draw  [->,bend left] (A) to  node [midway,above,align=center ]  {1} (C);
 \draw  [->,bend left] (C) to  node  [midway,above,align=center ] {1} (A);
 \draw  [->] (C) to  node [midway,above,align=center ]  {0} (D);

 \path (B) edge [loop above] node {0,1} (B);
 \path (D) edge [loop above] node {0,1} (D);

\end{tikzpicture}
\end{center}
\caption{Apwenian automaton} \label{fig:apw}
\end{figure}

In addition to the examples above, all ultimately periodic sequences are
$k$-automatic for all integers $k\geq 2$, see \cite[Theorem~5.4.2]{alsh2}.
Moreover, by  Cobham's theorem \cite[Theorem~11.2.1]{alsh2}, if a sequence $(s_n)$ is both $k$-automatic and $\ell$-automatic and $k$ and $\ell$ are multiplicatively 
independent,\footnote{Two integers $k$ and $\ell$ are {\em multiplicatively dependent} if $k^r=\ell^s$ for some positive integers $r$ and~$s$. Otherwise they are {\em multiplicatively independent}.} 
then $(s_n)$ is ultimately periodic.

For a prime power $k=q$, $k$-automatic sequences $(s_n)$ over the finite field\footnote{For a prime power $q$ we denote the finite field of size $q$ by $\F_q$.} $\Delta=\F_q$ can be characterized  
by a result of Christol, see \cite{christol} for prime $q$ and \cite{chka} for prime power $q$ as well as \cite[Theorem~12.2.5]{alsh2}.

\begin{theorem}\label{thm:christol}
Let 
\footnote{We denote by $\F_q \llbracket  x \rrbracket$ the ring of formal power series over $\F_q$.}
$$
G(x)=\sum_{n=0}^{\infty}s_nx^n \in \F_q \llbracket  x \rrbracket 
$$
be the {\em generating function} of the sequence $(s_n)$ over $\F_q$. 
Then $(s_n)$ is $q$-automatic if and only if $G(x)$ is algebraic over $\F_q(x)$, that is, there is a polynomial 
$h(x,y)\in \F_q[x,y]\setminus\{0\}$
such that $h(x,G(x))=0$. 
\end{theorem}
Note that for all $m=1,2,\ldots$ a sequence is $k$-automatic if and only if it is $k^m$-automatic by \cite[Theorem~6.6.4]{alsh2} and even a slightly more general version of Christol's result holds:
 For a prime $p$ and positive integers $m$ and $r$, $(s_n)$ is $p^m$-automatic over $\F_{p^r}$ if and only if $G(x)$ is algebraic over $\F_{p^r}(x)$.

\begin{example}
The generating function  $G(x)$ of the Thue-Morse sequence $(t_n)$ over $\F_2$ satisfies $h(x,G(x))=0$ with
\begin{equation}\label{eq:TH-equation}
h(x,y)=(x+1)^3y^2 + (x+1)^2y+x.
\end{equation}

The generating function  $G(x)$ of the Rudin-Shapiro sequence $(r_n)$ over $\F_2$ satisfies $h(x,G(x))=0$ with
\begin{equation}\label{eq:RS-equation}
h(x,y)=(x+1)^{5}y^2 + (x+1)^4 y + x^3.
\end{equation}

In general, for prime $p$ the generating function $G(x)$ of the $p$-ary pattern sequence $(p_n)$ over $\F_p$ with respect to the pattern $P$ of length $\ell$ satisfies $h(x,G(x))=0$ with
\begin{equation}\label{eq:pattern-equation}
h(x,y)=(x-1)^{p^\ell +p -1}y^p - (x-1)^{p^\ell} y - x^{a(P)}.
\end{equation}

The generating function $G(x)$ of the Rudin-Shapiro-like 
sequence $(\ell_n)$ over $\F_2$ defined by~\eqref{eq:rslike} satisfies $h(x,G(x))=0$ with
\begin{equation}\label{eq:rslh}
h(x,y)=(x+1)^8y^4+(x^6+x^5+x^2+x)y^2+(x+1)^4y+x^2,
\end{equation}
see \cite[Proof of Theorem 2]{suzeli}.

The generating function  $G(x)$ of the Baum-Sweet sequence $(b_n)$ over $\F_2$ satisfies \linebreak[4] $h(x,G(x))=0$ with
\begin{equation}\label{eq:BS-equation}
h(x,y)=y^3+xy+1.
\end{equation}

The generating function $G(x)$ of the characteristic sequence $(c_n)$ of sums of three 
squares~\eqref{eq:cndef} over $\F_2$ satisfies $h(x,G(x))=0$ with  
\begin{equation}\label{eq:cnh}
h(x,y)=(x+1)^8(y+y^4)+x^6+x^5+x^3+x^2+x,
\end{equation}
see \cite[Equation (7)]{howi}.

The generating function $G(x)$ of the regular paper-folding sequence $(v_n)$ over $\F_2$ satisfies $h(x,G(x))=0$ with
\begin{equation}\label{pfh}
h(x,y)=(x+1)^4(y^2+y)+x.
\end{equation}

The generating function $G(x)$ of the apwenian sequence $(w_n)$ over $\F_2$ defined by \eqref{apdef}
satisfies
\begin{equation}\label{aph} h(x,y)=(x+1)(xy^2+y)+1.
\end{equation}
\end{example}

\section{Linear complexity}\label{sec:linear_complexity}

The linear complexity is a figure of merit of pseudorandom sequences 
introduced to capture undesirable linear structure in a sequence. It originates in cryptography and provides a test  of randomness which is 
a standard tool to filter sequences with non-randomness properties and is implemented in many test suites such as NIST and TestU01~\cite{nist,testU01}.

\begin{definition}
The \emph{$N$th linear complexity} $L(s_n, N)$ of a sequence $(s_n)$ over $\F_q$ is the length~$L$ of a shortest linear recurrence relation satisfied by the first $N$ elements of $(s_n)$,
$$
 s_{n+L}=c_{L-1}s_{n+L-1}+\dots +c_1s_{n+1}+c_0s_n, \quad 0\leq n\leq N-L-1,
$$
for some $c_0,\ldots,c_{L-1}\in \F_q$.  
We use the convention that $L(s_n,N)=0$ if the first $N$ elements of $(s_n)$ are all zero and $L(s_n,N)=N$ if $s_0=\dots=s_{N-2}=0\ne s_{N-1}$.
The sequence~$(L(s_n,N))_{N=1}^\infty$
is called {\em linear complexity profile} of $(s_n)$ and  
$$
L(s_n)=\sup_{N\ge 1} L(s_n,N)
$$ 
is the {\em linear complexity} of $(s_n)$. 
\end{definition}

Clearly, $0\leq L(s_n,N) \leq N$ and $L(s_n,N)\leq L(s_n,N+1)$.

For truly random sequences $(s_n)$ the expected value of its $N$th linear complexity $L$ is 

$$\frac{N}{2}+O(1),
$$
see for example \cite[Theorem~10.4.42]{handbookFF}.
Deviations of order of magnitude $\log N$ must appear for infinitely many $N$.
More precisely, for a prime power $q$ consider the following probability measure of sequences over $\F_q$ determined by
\begin{equation}\label{eq:prob_space}
\mathbb{P}\left[(s_n)\in \F_q^{\infty}: (s_0,\dots, s_{\ell-1})=  (c_0,\dots, c_{\ell-1}) \right]=q^{-\ell},  \quad     c_0,\dots, c_{\ell-1}\in \F_q.
\end{equation}
Then we have the following result on the deviation from the expected value, see\cite[Theorem~10]{Niederreiter88}.
\begin{theorem}
We have
$$
\limsup_{N\rightarrow \infty }\frac{L(s_n, N) -N/2}{\log N}=\frac{1}{2 \log q}, 
$$
and 
$$
\liminf_{N\rightarrow \infty }\frac{L(s_n, N) -N/2}{\log N}=\frac{-1}{2 \log q}
$$
 with probability one with respect to the probability measure \eqref{eq:prob_space}.
\end{theorem}

It is well-known \cite[Lemma~1]{Niederreiter88-b} that $L(s_n)<\infty$ if and only if $(s_n)$ is ultimately periodic, that is,
its generating function is rational: $G(x)=g(x)/f(x)$ with polynomials $g(x),f(x)\in \F_q[x]$.

The $N$th linear complexity is a measure for the unpredictability of a sequence. 
A large $N$th linear complexity, up to sufficiently large $N$,
is necessary, but not sufficient, for cryptographic applications. Sequences of small linear complexity are also weak in view of Monte-Carlo methods, see \cite{do,domewi,dowi,DorferWinterhof}. 
For more background on linear complexity and related measures of pseudorandomness we refer to 
\cite[Section~10.4]{handbookFF}
and \cite{Niederreiter2003,towi,wi1}.

M\'erai and Winterhof \cite{mewi18} showed that automatic sequences which are not ultimately periodic possess large $N$th linear complexity.

\begin{theorem}\label{thm:MW-lin-compl-gen}
 Let $q$ be a prime power and $(s_n)$ be a $q$-automatic sequence over $\F_q$ which is not ultimately periodic. 
Let $h(x,y)=h_0(x)+h_1(x)y+\dots + h_d(x)y^d\in\F_q[x,y]$ be a non-zero polynomial with $h(G(x),x)=0$ with no rational zero.
Put $$
M=\max_{0\leq i\leq d}\{\deg h_i-i\}.
$$
Then we have
\[
 \frac{\displaystyle N-M}{d}\leq L(s_n,N)\leq \frac{\displaystyle (d-1)N+M+1}{d}.
\]
\end{theorem}
See also \cite{XingLam} for the special case $d=2$.

The idea of the proof of Theorem \ref{thm:MW-lin-compl-gen} is that small $N$th linear complexity profile gives a good rational approximation to the generating function. However, transcendental elements over $\F_q(x)$ are not well-approximated.

Namely, since $(s_n)$ is not ultimately periodic, $G(x)=\sum_{n=0}^\infty s_nx^n\not\in \F_q(x)$ is not rational by \cite[Lemma~1]{Niederreiter88-b}.

Let $g(x)/f(x)\in\F_q(x)$ be a rational zero of $h(x,y)$ modulo $x^N$ with $\deg(f)\le L(s_n,N)$ and $\deg(g)<L(s_n,N)$. 
More precisely, put $L=L(s_n,N)$. Then we have  
$$
\sum_{\ell=0}^Lc_{\ell}s_{n+\ell}=0 \quad \mbox{for }0\le n\le N-L-1
$$ 
for some $c_0,\ldots,c_{L}\in \F_p$ with $c_L=-1$.
Take
\[
f(x)=\sum_{\ell=0}^L c_{\ell}x^{L-\ell} 
\]
and 
\[
g(x)=\sum_{m=0}^{L-1}\left(\sum_{\ell=L-m}^Lc_\ell s_{m+\ell-L}\right)x^m 
\]
and verify 
$$
f(x)G(x)\equiv g(x)\bmod x^N.
$$
Then
\[
 h_0(x)f^d(x)+h_1(x)g(x)f^{d-1}(x)+\dots + h_d(x)g(x)^d=K(x) x^N.
\]
Here $K(x)\neq 0$ since $h(x,y)$ has no rational zero. Comparing the degrees of both sides we get
\[
 dL+M\geq N
\]
which gives the lower bound.

The upper bound for $N=1$ is trivial. For $N\geq 2$ the result follows from the well-known bound, see for example \cite[Lemma 3]{DorferWinterhof},
\[
 L(s_n,N)\leq \max\left\{L(s_n,N-1), N-L(s_n,N-1) \right\}
\]
by induction.

The bound in Theorem \ref{thm:MW-lin-compl-gen} combined with \eqref{eq:TH-equation}-\eqref{aph}
gives the following estimates for the $N$th linear complexity of the Thue-Morse sequence $(t_n)$ defined by~\eqref{tmdef}
\begin{equation}\label{eq:lin-compl_TM}
\left\lceil \frac{N-1}{2} \right\rceil\leq  L(t_n,N)\leq \left\lfloor \frac{N}{2} \right\rfloor+1,
\end{equation}
of the Rudin-Shapiro sequence $(r_n)$ defined by \eqref{rsdef} and the regular paper-folding sequence $(v_n)$ defined by \eqref{pfdef}
\begin{equation}\label{eq:lin-compl_RS}
\left\lceil\frac{N-3}{2}\right\rceil\leq L(r_n,N), L(v_n,N)\leq \left\lfloor\frac{N}{2}\right\rfloor+2,
\end{equation}
of the $p$-ary pattern sequence $(p_n)$ defined by \eqref{eq:recurrence} with any pattern $P$ of length $\ell$
$$
\left\lceil\frac{N+1}{p}\right\rceil-p^{\ell-1}\leq L(p_n,N)\leq \left\lfloor\frac{(p-1)N}{p}\right\rfloor+p^{\ell-1},
$$
of the Rudin-Shapiro-like sequence $(\ell_n)$ defined by \eqref{eq:rslike}
$$
\left\lceil \frac{N}{4}\right\rceil-1\le L(\ell_n,N)\le \left\lfloor\frac{3N+5}{4}\right\rfloor,
$$
of the Baum-Sweet sequence $(b_n)$ defined by \eqref{bsdef}
$$
\left\lceil\frac{N}{3}\right\rceil\leq L(b_n,N)\leq \left\lfloor\frac{2N+1}{3}\right\rfloor, 
$$
of the characteristic sequence $(c_n)$ of sums of three squares defined by \eqref{eq:cndef}
$$
\left\lceil\frac{N-7}{4}\right\rceil\le L(c_n,N)\le \left\lfloor\frac{3N}{4}\right\rfloor+2
$$
and of the apwenian sequence $(w_n)$ defined by \eqref{apdef}
\begin{equation}\label{aplin} L(w_n,N)=\left\lfloor \frac{N+1}{2}\right\rfloor.
\end{equation}
Note that the bound \eqref{eq:lin-compl_TM} is also true for the dual $(t'_n)$ of the Thue-Morse sequence, that is, $t'_n=1-t_n$, and apwenian sequences are characterized by the property \eqref{aplin}, see \cite{alhani}. Note that not all apwenian sequences are automatic.

The bounds \eqref{eq:lin-compl_TM} for the Thue-Morse sequence and \eqref{eq:lin-compl_RS}
for the Rudin-Shapiro sequence are optimal.
Using the continued fraction expansions of their generating functions, M\'erai and Winterhof
\cite{mewi18} determined the exact value of the $N$th linear complexity profiles of the Thue-Morse and Rudin Shapiro sequence.
\begin{theorem}
The $N$th linear complexity of the Thue-Morse sequence is
$$
L(t_n,N)=2\left \lfloor \frac{N+2}{4}\right\rfloor, \quad N=1,2,\dots
$$
and the $N$th linear complexity of the Rudin-Shapiro sequence is 
$$L(r_n,N)=\left\{\begin{array}{cc} 6\left\lfloor N/12\right\rfloor+4, & N\equiv 4,5,6,7,8,9\bmod 12,\\
6\left\lfloor (N+2)/12\right\rfloor, & \mbox{otherwise}.
\end{array}\right.$$
\end{theorem}
The result can be extended to binary pattern sequences $(p_n)$ defined by \eqref{eq:recurrence} with the all one pattern of length $\ell\ge 3$, that is, $a=2^\ell-1$.

It follows from Theorem~\ref{thm:MW-lin-compl-gen}, that if an automatic sequence is not ultimately periodic and its generating function has a quadratic minimal polynomial, that is $d=2$ in Theorem~\ref{thm:MW-lin-compl-gen}, then the deviation of the $N$th linear complexity from its expected value $N/2$ is bounded by $(M+1)/2$,
$$
\left|L(s_n,N)-\frac{N}{2}\right|\leq \frac{M+1}{2}.
$$
Such sequences 
are said to have \emph{almost perfect} or \emph{$(M+1)$-perfect} linear complexity profile, 
see \cite{Niederreiter88-b,alhani}.

Apwenian sequences are those sequences having $1$-perfect or just {\em perfect} linear complexity profile. The bounds \eqref{eq:lin-compl_TM} and
\eqref{eq:lin-compl_RS} imply that the Thue-Morse sequence has $2$-perfect linear complexity profile and the Rudin-Shapiro sequence and the paper-folding sequence both have $4$-perfect linear complexity profile.

\begin{figure}[ht]
\begin{center}
\includegraphics[scale=.51]{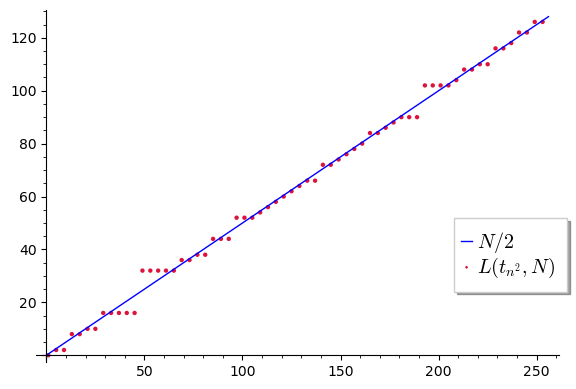} 
\includegraphics[scale=.51]{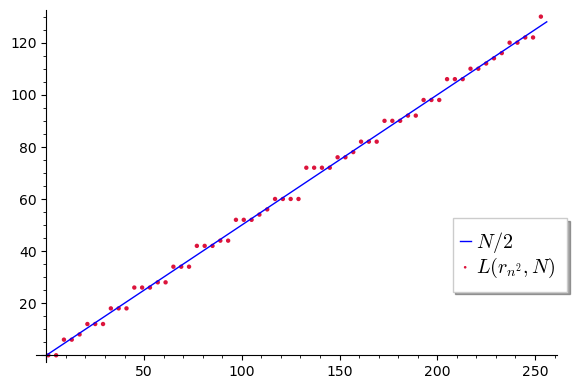} 
\end{center}
\caption{The $N$th linear complexity of the Thue-Morse (left) and 
Rudin-Shapiro (right) sequence  along squares.}
\label{fig:L}
\end{figure}

Although automatic sequences have some good pseudorandom properties including a desirable linear complexity profile, these sequences have also some strong non-randomness properties, see Sections~\ref{sec:correlation}, \ref{sec:expansion_complexity} and \ref{sec:normality} below.
Such randomness flaws may be avoided considering subsequences of automatic sequences. For example, the Thue-Morse and Rudin-Shapiro sequences along squares are not automatic, see Section~\ref{sec:normality} below, and seem to have $N$th linear complexity $N/2+O(\log N)$, see Figure~\ref{fig:L}. 

\begin{problem}
Prove that the $N$th linear complexities of the Thue-Morse and Rudin-Shapiro sequences along squares satisfy
\footnote{$f(k)=o(g(k))$ is equivalent to $f(k)/ g(k)\rightarrow 0$ as $k\rightarrow \infty$.}
$$
L(t_{n^2},N)=\frac{N}{2}+o(N)\quad \mbox{and}\quad L(r_{n^2},N)=\frac{N}{2}+o(N).
$$
\end{problem}
We remark, that lower bounds on the $N$th linear complexities of $(t_{n^2})$ and $(r_{n^2})$ of order of magnitude $\sqrt{N}$ follow from Theorem~\ref{thm:suwi} and \eqref{eq:ML} in the next section.

Additional to these examples, the same problem is also open for other subsequences such as along other polynomial values, along primes etc.

\section{Maximum order complexity}\label{sec:max-order_complexity}

Maximum order (or nonlinear) complexity is a refinement of the linear complexity considering not only linear but \emph{any} recurrence relation.

\begin{definition}
The {\em $N$th maximum order complexity}  $M(s_n,N)$ is the smallest positive integer~$M$ with 
 $$
 s_{n+M}=f(s_{n+M-1},\ldots,s_n),\quad 0\le n\le N-M-1,
 $$
 for some mapping $f:\F_2^M \rightarrow \F_2$. 
 The sequence $(M(s_n,N))_{N=1}^\infty$ is called {\em maximum order complexity profile}.
\end{definition}

Obviously, we have 
\begin{equation}\label{eq:ML}
 M(s_n,N)\le L(s_n,N)
\end{equation} 
and the maximum order complexity is a finer measure for the unpredictability of a sequence than the linear complexity.
However, often the linear complexity is easier to analyze both theoretically and algorithmically. 

Clearly, a sufficiently large maximum order complexity is needed for unpredictability and suitability in cryptography. 
However, sequences of very large maximum order complexity have also a very large autocorrelation or correlation measure of order $2$, see \eqref{C2M} below, and are not suitable for many applications including cryptography, radar, sonar and wireless communications.

The maximum order complexity was introduced by Jansen in \cite[Chapter 3]{ja}, see also~\cite{jabo}.
The typical value for the $N$th maximum order complexity is of order of magnitude $\log N$, see
\cite{ja,jabo}. An algorithm for calculating the maximum order complexity profile of linear time and memory was presented by 
Jansen \cite{ja,jabo} using the graph algorithm introduced by Blumer et al.\ \cite{bl}.

The maximum order complexity of the Thue-Morse sequence was determined in \cite[Theorem 1]{suwi19}.
\begin{theorem}\label{thm:suwi-tm}
  For $N\ge 4$,
  the $N$th maximum order complexity of the Thue-Morse sequence~$(t_n)$
  satisfies
$$M(t_n,N)=2^\ell+1,$$
where
$$\ell=\left\lceil \frac{\log (N/5)}{\log 2}\right\rceil.$$
\end{theorem}
 It is easy to see that
\begin{equation}\label{tmmoc}\frac{N}{5}+1\le M(t_n,N)\le 2\frac{N-1}{5}+1
 \quad\text{for}\;\; N\ge 4.
\end{equation}
In Section~\ref{sec:correlation} we will see that such a large maximum order complexity points to undesirable structure in a sequence. 

The $N$th maximum order complexity of the Rudin-Shapiro sequence and some generalizations is also of order of magnitude $N$, see \cite[Theorem 2]{suwi19}. In particular we have
\begin{equation}\label{rsmoc}
M(r_n,N)\ge \frac{N}{6}+1,\quad N\ge 4.
\end{equation}

The maximum order complexity of the subsequences of the Thue-Morse and the Rudin-Shapiro sequence along squares are still
large enough, see \cite{suwi}.
\begin{theorem}\label{thm:suwi}
 The $N$th maximum order complexities $M(t_{n^2},N)$ and $M(r_{n^2},N)$
 of the subsequences $(t_{n^2})$ and $(r_{n^2})$ of the Thue-Morse and the Rudin-Shapiro sequence along squares satisfy
\begin{align*}
   &M(t_{n^2},N)\geq \sqrt{\frac{2N}{5}},\quad N\ge 21,   \quad  \mbox{and} \\
  &M(r_{n^2},N)\geq \sqrt{\frac{N}{8}},\quad N\ge 64.
\end{align*}   
\end{theorem}
We sketch the proof.
First, let $t$ be the length of the longest subsequence of~$(t_{n^2})$ that
occurs at least twice with different successors among the first $N$ sequence elements. Then 
$M(t_{n^2},N)\ge t + 1$. Hence
 the first inequality follows from
$$
t_{(i+2^{\ell+1})^2}=t_{(i+2^{\ell+2})^2},\quad i=0,1,\ldots,\left\lfloor \sqrt{2^{\ell+2}-1}\right\rfloor$$
$$\mbox{and}\quad t_{(2^{\ell}+2^{\ell+1})^2}\ne t_{(2^\ell+2^{\ell+2})^2},
$$
 which can be shown by induction over $\ell\ge 2$, 
where $\ell$ is defined by $5\cdot 2^\ell<N\le 5\cdot 2^{\ell+1}$.
 
 The second bound follows from
 $$r_{(i+2^{\ell+3})^2}=r_{(i+2^{\ell+4})^2},\quad i=0,1,\ldots,\left\lfloor\sqrt{2^{\ell+3}-1}\right\rfloor,$$
 $$\mbox{and}\quad
 r_{(2^{\ell+2}+2^{\ell+3})^2}\ne r_{(2^{\ell+2}+2^{\ell+4})^2},$$
 where $\ell$ is defined by $2^{\ell+5}\le N<2^{\ell+6}$.

\begin{figure}[ht]
\begin{center}
\includegraphics[scale=.5]{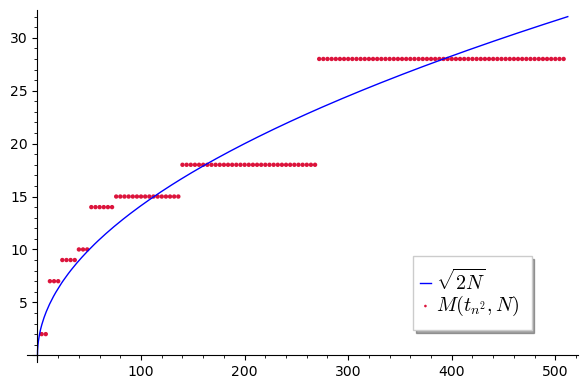}
\includegraphics[scale=.5]{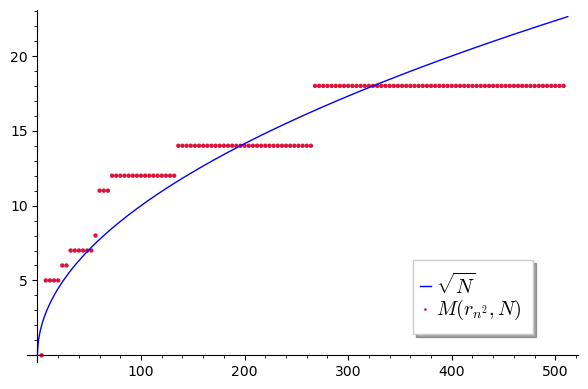}
\end{center}
\caption{The $N$th maximum order complexity of the Thue-Morse (left) and Rudin-Shapiro (right) sequence along squares.}
\label{maxordersquare}
\end{figure}

Figure \ref{maxordersquare} suggests that $\sqrt{N}$ is the right order of magnitude for the $N$th maximum order complexities of $(t_{n^2})$ and $(r_{n^2})$.
For $N\ge 2^{2\ell+2}$ the same lower bound $\sqrt{N/8}$ is true for binary pattern sequences along squares with the all one pattern of length $\ell$, that is, $a=2^\ell-1$ for~$\ell\ge 3$, see \cite{suwi}.

This result was extended by Popoli \cite{po} to sequences along polynomial values of higher degrees $d$. However, the lower bounds are of order of magnitude $N^{1/d}$.
Note that no better lower bounds are known for the $N$th linear complexity of these subsequences of automatic sequences.

The problem for other subsequences is still open as for example for subsequences along primes.
\begin{problem}
 Study the maximum-order complexity of the subsequences of the Thue-Morse and the Rudin-Shapiro sequence along primes.
\end{problem}

The maximum order complexity of other automatic sequences has also been studied. 
Sun, Zeng and Lin \cite{suzeli} showed that the $N$th maximum order complexity of the Rudin-Shapiro-like sequence $(\ell_n)$ defined by \eqref{eq:rslike} is of order of magnitude~$N$.

\bigskip

We remark, that in addition to automatic sequences based on the $k$-ary expansion~\eqref{eq:k-ary} of integers, one can consider analogously sequences using other numeration systems.

In particular, consider the \emph{Fibonacci numbers} defined by
$$
F_0=0,~F_1=1\quad \mbox{and}\quad F_n=F_{n-1}+F_{n-2}\mbox{ for }n\ge 2.
$$
Then the unique, see for example \cite[Theorem~3.8.1]{alsh2}, {\em Zeckendorf expansion} or {\em Fibonacci expansion}, of a positive integer~$n$ is
$$
n=\sum_{i=0}^\infty e_i F_{i+2},\quad \mbox{where }e_i\in \{0,1\} \mbox{ and } e_ie_{i+1}=0\mbox{ for }i=0,1,\ldots
$$
Analogously to the Thue-Morse sum-of-digits sequence $(t_n)$ and the Rudin-Shapiro sequence~$(r_n)$
which can be defined by~\eqref{sumofdigitsdef} we can define and study the {\em Zeckendorf 
sum-of-digits sequences modulo $2$} $(z_n)$ and $(u_n)$ defined by
\begin{equation}\label{zeck}
z_n=\sum_{i=0}^\infty e_i \bmod 2
\quad\mbox{and}\quad
u_n=\sum_{i=0}^\infty e_ie_{i+2}\bmod 2.
\end{equation}
Very recently the maximum-order complexity of $(z_n)$ and its subsequences along polynomial values has been studied by Jamet, Popoli and Stoll in \cite{japost}. 
A lower bound on $M(u_n,N)$ and some generalizations can be obtained along the same lines and will be contained in Popoli's thesis. 

\section{Well-distribution and correlation measures}\label{sec:correlation}

Mauduit and S\'ark\"ozy \cite{masa} introduced two measures
of pseudorandomness for finite sequences over $\{-1,+1\}$, the well-distribution measure and the correlation measure of order $k$. We adjust these definitions to infinite binary sequences $(s_n)$ over $\F_2$.

\begin{definition}
The {\em $N$th well-distribution measure} of $(s_n)$ is defined as
$$W(s_n,N)=\max_{a,b,t}\left|\sum_{j=0}^{t-1} (-1)^{s_{a+jb}}\right|,$$
where the maximum is taken over all integers $a,b,t$, $b\ge 1$, such that $0\le a \le a+(t-1)b\le N-1$.
\end{definition}

The well-distribution measure provides information on the balance,
\footnote{Note that the term {\em balanced} is used with a different meaning in  combinatorics on words, see for example \cite[Definition~10.5.4]{alsh2}.}
that is the distribution of zeros and ones, along arithmetic progressions.
For random sequences it is expected to be small. More precisely, Alon et al.\ \cite[Theorem~1]{alko} proved the following result on the typical value of the well-distribution measure.

\begin{theorem}
For all $\varepsilon>0$, there are numbers $N_0=N_0(\varepsilon)$ and $\delta=\delta(\varepsilon)>0$ such that for $N\ge N_0$ we have 
$$\delta \sqrt{N}<W(s_n,N) <\frac{\sqrt{N}}{\delta}$$
with probability at least $1-\varepsilon$ with respect to the probability measure \eqref{eq:prob_space}.
\end{theorem}

Moreover, Aistleitner \cite{aistleitner} showed that
there exists a continuous limit distribution 
of~$\frac{W(s_n,N)}{\sqrt{N}}$.
More precisely, 
for any $t\in \R$ the limit
$$
F(t)=\lim_{N\rightarrow \infty} {\mathbb P}\left(\frac{W(s_n,N)}{\sqrt{N}}\le t\right) 
$$
exists and satisfies 
$$
\lim_{t\rightarrow \infty}t(1-F(t))e^{t^2/2}=\frac{8}{\sqrt{2\pi}},
$$
with respect to the probability measure~\eqref{eq:prob_space}.

\begin{definition}
For $k\geq 1$,
the {\em $N$th correlation measure of order~$k$} of a binary sequence~$(s_n)$ is 
$$C_k(s_n,N)=\max_{M,D}\left|\sum^{M-1}_{n=0}(-1)^{s_{n+d_1}}\cdots (-1)^{s_{n+d_k}}\right|,$$ 
where the maximum is taken over all $D=(d_1,d_2,\ldots,d_k)$ with integers satisfying
$0\le d_1<d_2<\cdots<d_k$ and $1\le M\le N-d_k$.
\end{definition}

The correlation measure of order $k$ provides information about the similarity of parts of the sequence and their shifts.
For a random sequence this similarity and thus the correlation measure of order $k$ is expected to be small.
More precisely,  Alon et al.\ \cite[Theorem 2]{alko} proved the following result on the typical
value of the correlation measure of 
order $k$.
\begin{theorem}
For any $\varepsilon>0$, there exist an $N_0=N_0(\varepsilon)$ such that for all $N\ge N_0$
we have for a randomly chosen sequence $(s_n)$ and any $k$ with $2\le k\le N/4$,
$$
\frac{2}{5}\sqrt{N\log{N\choose k}}<C_k(s_n,N)<\frac{7}{4}\sqrt{N\log{N\choose k}}
$$
with probability at least $1-\varepsilon$ with respect to the probability measure \eqref{eq:prob_space}.
\end{theorem}

Moreover, Schmidt \cite[Theorem 1.1]{schmidt} showed, that for fixed $k$, we have
$$
\lim_{N\rightarrow \infty}\frac{C_k(s_n,N)}{\sqrt{2N\log \binom{N}{k-1}}} =1
$$
with probability $1$ with respect to the probability measure \eqref{eq:prob_space}.

A large well-distribution measure implies a large correlation measure of order $2$. More precisely we have by \cite[Theorem 1]{ms03}
\footnote{$f(k)=O(g(k))$ is equivalent to $|f(k)|\le c g(k)$ for some constant $c>0$.}
$$
W(s_n,N)=O\left(\sqrt{NC_2(s_n,N)}\right).
$$

Mauduit and Sárközy \cite{masa98} obtained bounds on the well-distribution measure and correlation measure of order $2$ of Thue-Morse sequence $(t_n)$ and Rudin-Shapiro sequence~$(r_n)$. 

For example, as a consequence of the bound
\begin{equation}\label{eq:exp-TM}
\left|\sum_{n=0}^{N-1}(-1)^{t_n}z^n\right|\leq (1+\sqrt{3}) N^{\log 3/\log 4},\quad |z|=1,
\end{equation}
of Gel'fond~\cite[p.~262]{gelfond}, see \cite{foma} for the explicit constant $1+\sqrt{3}$, they obtained a bound on~$W(t_n,N)$. 
\begin{theorem}\label{thm:masa-thm1} 
 We have 
 $$
 W(t_n,N) \leq 2(1+\sqrt{3}) N^{\log 3/\log 4}.
 $$
\end{theorem}
Also, using the bound 
\begin{equation}\label{eq:exp-RS}
\left|\sum_{n=0}^{N-1}(-1)^{r_n}z^n\right|\leq (2+\sqrt{2}) N^{1/2},\quad |z|=1,
\end{equation}
obtained by Rudin~\cite{rudin} and Shapiro~\cite{shapiro}, see also \cite[Theorem 3.3.2]{alsh2},
they proved a bound on~$W(r_n,N)$.
\begin{theorem}\label{thm:masa-thm2} 
 We have 
 $$
 W(r_n,N)\leq 2(2+\sqrt{2}) N^{1/2}.
 $$
\end{theorem}
In general, following the proofs of \cite{masa98} we get
\begin{equation}\label{eq:W-idea}
W(s_n,N)
=O\left(
 \sup_{|z|=1, m\le N}\left|\sum_{n=0}^{m-1} (-1)^{s_n}z^n\right| \right)
\end{equation}
and thus Theorems \ref{thm:masa-thm1} and \ref{thm:masa-thm2} follow, up to the constant, from \eqref{eq:exp-TM} and \eqref{eq:exp-RS}.

However, for $(t_n)$ and $(r_n)$ Mauduit and Sárközy \cite{masa98} detected non-randomness properties by showing that the correlation measure of order $2$ of these sequences is large.
\begin{theorem}
 We have 
 \begin{equation}\label{ctm} C_2(t_n,N)>\frac{N}{12},\quad N\ge 5,
 \end{equation}
 and
 \begin{equation}\label{crs}
  C_2(r_n,N)>\frac{N}{6},\quad N\ge 4.
 \end{equation}
\end{theorem}

Mérai and Winterhof~\cite{mewi2} showed that all automatic 
sequences share the property of having a large correlation measure of order $2$. They provided the following lower bound in terms of the defining automaton.
\begin{theorem}\label{thm:gen_lowe_bound_on_C}
Let $(s_n)$ be a $k$-automatic binary sequence generated by the finite automaton $(Q,\Sigma,\delta,q_0,\varphi,\{0,1\})$.
Then
 \[
  C_2(s_n,N)\geq \frac{N}{k(|Q|+1)} \quad \text{for } N\geq k(|Q|+1).
 \]
\end{theorem}

This result applied to $(t_n)$ and $(r_n)$ gives the following bounds
$$
C_2(t_n,N)\geq \frac{N}{6}, \quad N\geq 6, \quad\text{and}\quad  C_2(r_n,N)\geq \frac{N}{10}, \quad N\geq 10,
$$
which improves \eqref{ctm}.

\begin{figure}[ht]
\begin{center}
\includegraphics[scale=.51]{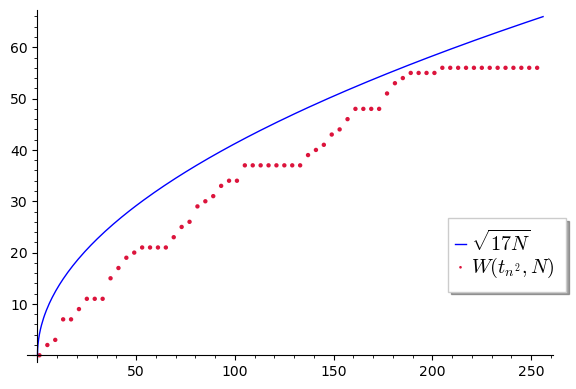} 
\includegraphics[scale=.51]{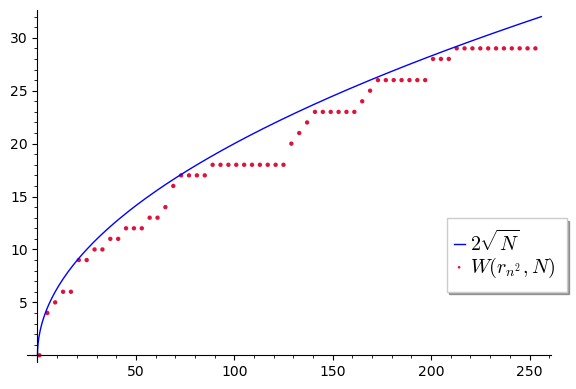} 
\end{center}
\caption{The $N$th well-distribution measure of the Thue-Morse (left) and Rudin-Shapiro (right) sequence along squares.}
\label{fig:well}
\end{figure}

\begin{figure}[ht]
\begin{center}
\includegraphics[scale=.51]{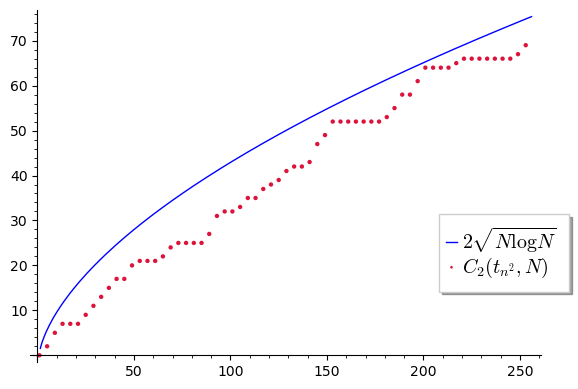}
\includegraphics[scale=.51]{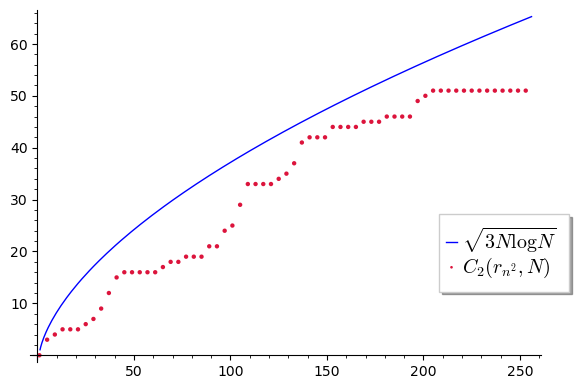} 
\end{center}
\caption{The $N$th second order correlation measure of the Thue-Morse (left) and Rudin-Shapiro (right) sequence along squares.}
\label{fig:cor}
\end{figure}

Figures~\ref{fig:well} and \ref{fig:cor} may lead to the conjecture that well-distribution measure and correlation measure of order $2$ of both $(t_{n^2})$ and $(r_{n^2})$ are of order of magnitude $N^{1/2}$ and $(N\log N)^{1/2}$, respectively.

\begin{problem}
For fixed $k=2,3,\ldots$
show that 
$$
C_k(t_{n^2},N)=o(N)\quad \mbox{and}\quad C_k(r_{n^2},N)=o(N).
$$
\end{problem}

Mauduit and Rivat \cite{mari18} showed that  
$$
\left|\sum_{n=0}^{N-1}(-1)^{r_{n^2}}z^n\right|= O\left( N^{1-\eta}\right),\quad |z|=1,\quad \mbox{for some  $\eta>0$},
$$
which, together with \eqref{eq:W-idea}, gives
a bound on $W(r_{n^2},N)$ of the same order of magnitude. 
More precisely, \cite{mari18} deals with the more general case of binary pattern sequences $(p_n)$ defined by \eqref{eq:pattern-equation} with either the all one pattern of length $k\ge 2$, that is, $a=2^k-1$, or the patterns $10\ldots 01$ of length $k\ge 3$, that is, $a=2^{k-1}+1$, and the constants depend on $k$.
For the Thue-Morse sequence along squares $(t_{n^2})$ one
can easily derive a nontrivial bound on
$$
\left|\sum_{n=0}^{N-1}(-1)^{t_{n^2}}z^n\right|,\quad |z|=1,
$$ 
and thus on $W(t_{n^2},N)$ since
the proof of \cite[Th\'eor\`eme 1]{mari1} for $z=1$
works also for $z\neq 1$ since after applying a variant of the van der Corput inequality, \cite[Lemma 15]{mari1}, we get an expression which does not depend on the variable $z$ anymore, that is, the same expression as for $z=1$.

Theorem~\ref{thm:suwi-tm} in Section~\ref{sec:max-order_complexity} above shows that the Thue-Morse sequence has maximum order complexity $M(t_n,N)$ of order of magnitude $N$.
Although a large maximum order complexity is desired it should be not too large since otherwise the correlation measure of order $2$ is large.
Namely, we have 
   \begin{equation}\label{C2M} C_2(s_n,N)\ge M(s_n,N)-1
   \end{equation}
   since by \cite[Proposition 3.1]{ja} there exist $0\le n_1<n_2\le N-M(s_n,N)-1$ with
   $$s_{n_1+i}=s_{n_2+i},\quad i=0,\ldots,M(s_n,N)-2,\quad \mbox{but }s_{n_1+M(s_n,N)-1}\ne s_{n_2+M(s_n,N)-1}$$
   and thus
   $$M(s_n,N)-1=\sum_{i=0}^{M(s_n,N)-2}(-1)^{s_{n_1+i}+s_{n_2+i}}\le C_2(s_n,N).$$
Combining \eqref{tmmoc} and \eqref{C2M} we get for the Thue-Morse sequence
$$
C_2(t_n,N)\ge \frac{N}{5},\quad N\ge 4,
$$
which further improves the constant in \eqref{ctm}.
Combining \eqref{rsmoc} and \eqref{C2M} recovers~\eqref{crs}.
The correlation measure of order $2$ with bounded lags of some generalizations of the Rudin-Shapiro sequence has recently been studied in \cite{mastta}.

In contrast to the Thue-Morse und Rudin-Shapiro sequence, the well-distribution measure of some other binary automatic sequences is very large.
For example,  
the Baum-Sweet sequence $(b_n)$, 
the characteristic sequence $(c_n)$ of the sums of three squares, 
the paper-folding sequence $(v_n)$ 
and the apwenian sequence $(w_n)$ defined by \eqref{apdef} are very unbalanced and thus
have all well-distribution measure of order of magnitude $N$.
However, it seems to be interesting to study the well-distribution measure for
arbitrary apwenian sequences.
For the Rudin-Shapiro like sequence $(\ell_n)$ defined by \eqref{eq:rslike} Lafrance, Rampersad and Yee \cite{laraye} proved
$$\liminf_{N\rightarrow \infty} \frac{\sum_{n=0}^{N-1}(-1)^{\ell_n}}{\sqrt{N}}=\frac{\sqrt{3}}{3}
\quad \mbox{and}\quad
\limsup_{N\rightarrow \infty} \frac{\sum_{n=0}^{N-1}(-1)^{\ell_n}}{\sqrt{N}}=\sqrt{2}.
$$
However, a bound on $W(\ell_n,N)$ is not known and in contrast to \eqref{eq:exp-RS} for the Rudin-Shapiro sequence~$(r_n)$, for $(\ell_n)$
the absolute values
$$\left|\sum_{n=0}^{N-1}(-1)^{\ell_n} z^n\right|$$
can be of much larger order of magnitude than $\sqrt{N}$ for some $z$ with $|z|=1$, see \cite[Theorem~2]{al16} as well as \cite{chgr}.

Finally, we remark that the result of Theorem~\ref{thm:gen_lowe_bound_on_C} provides an estimate on the \emph{state complexity} of sequences in terms of the correlation measure of order $2$.

\begin{definition}
Let $k\geq 2$. Then the \emph{$N$th state complexity} $SC_k(s_n,N)$ of a sequence~$(s_n)$ over $\F_2$ is the minimum of the number of states of any finite $k$-automaton which generates the first $N$ sequence elements.
\end{definition}

\begin{cor}
 Let $(s_n)$ be a binary sequence. Then for all $k\geq 2$ we have
 \[
  SC_k(s_n,N)\geq\frac{N}{k\cdot C_2(s_n,N)}-1 \quad \text{for } N\geq 3.
 \]
\end{cor}

\section{Expansion complexity}\label{sec:expansion_complexity}

Theorem~\ref{thm:MW-lin-compl-gen} indicates that automatic sequences possess good properties in terms of the linear complexity profile. 
However, the results of Section~\ref{sec:correlation} show that these sequences have a serious lack of 
pseudorandomness.
Diem~\cite{di} showed that these sequences are 
not just statistically auto-correlated, but 
are completely predictable from a relatively short initial segment. 
He introduced the notion of \emph{expansion complexity} 
to turn such security flaw into a quantitative form.

\begin{definition}
 Let $(s_n)$ be a sequence over $\F_q$ with generating function
 $$
 G(x)=\sum_{n=0}^{\infty}s_nx^n \in \F_q \llbracket  x \rrbracket .  
 $$
 For a positive integer $N$, the \emph{$N$th expansion complexity $E(s_n,N)$ of $(s_n)$} is $E(s_n,N)=0$ if $s_0=\dots=s_{N-1}=0$ and otherwise the least total degree of a non-zero polynomial $h(x,y)\in\F_q[x,y]$ such that
 \begin{equation}\label{eq:h}
 h(x,G(x))\equiv 0 \bmod x^N.
 \end{equation}
 The sequence $(E(s_n,N))_{N=1}^\infty$ is called \emph{expansion complexity profile of $(s_n)$} and
 $$
 E(s_n)=\sup_{N\geq 1} E(s_n,N)
 $$
 is the \emph{expansion complexity of $(s_n)$}.
\end{definition}

By Christol's Theorem~\ref{thm:christol}, a sequence is automatic if and only if its expansion complexity is finite. For example, we have for the Thue-Morse sequence $(t_n)$, the Rudin-Shapiro sequence $(r_n)$, the $p$-ary pattern sequence $(p_n)$, the Baum-Sweet sequence $(b_n)$, the Rudin-Shapiro like sequence $(\ell_n)$
and the characteristic sequence $(c_n)$ of sums of three squares that
$$
E(t_n)= 5, \quad E(r_n)= 7, \quad E(p_n)\le p^\ell+2p-1, \quad E(b_n)= 3, \quad E(\ell_n)\le 12
$$
$$
E(c_n)\le 12, \quad E(v_n)= 6 \quad \mbox{and}\quad E(w_n)= 4,
$$
which follows from \eqref{eq:TH-equation}, \eqref{eq:RS-equation}, \eqref{eq:pattern-equation},  \eqref{eq:rslh},
\eqref{eq:BS-equation}, \eqref{eq:cnh}, \eqref{pfh} and \eqref{aph}. The equalities follow from
the fact that there is no lower degree polynomial with such property since~$h(x,y)$ is irreducible in these cases, see \cite[Proposition~4]{di}. 

Diem showed \cite{di} that if a sequence has small expansion complexity, then 
long parts of such sequences can be computed efficiently from short ones. We summarize his 
results.
\begin{theorem}\label{thm:Diem}
 Let $(s_n)$ be a sequence over $\F_q$ with expansion complexity $E(s_n)=d$. From the first $d^2$ elements, one can compute an irreducible polynomial $h(x,y)\in \F_q[x,y]$ of degree $\deg h \leq d$ with $h(x,G(x))=0$ in polynomial time in $d\cdot \log q$.
 
 Moreover, an initial segment of the sequence of length $M>N$ can be determined from $h$ and the $d^2$ initial values in polynomial time in $d \cdot \log q$ and in linear time in~$M$.
\end{theorem}

Theorem~\ref{thm:Diem} shows that automatic sequences have a strong non-randomness property. The expansion complexity profile is defined to capture such non-randomness property locally, that is for initial segments of sequences. 

For the $N$th expansion complexity, we have the trivial bound $E(s_n,N)\leq N-1$ realized by the polynomial
$$
h(x,y)=y-\sum_{n=0}^{N-1}s_nx^n.
$$
Moreover, one can show the stronger upper bound
\begin{equation}\label{eq:exp_upper}
\binom{E(s_n,N)+1}{2}\leq N,
\end{equation}
which holds for all sequence $(s_n)$ and all $N\geq 1$, see \cite[Theorem~1]{GoMeNi}.

The $N$th expansion complexity of random sequences is concentrated to its upper bound \eqref{eq:exp_upper}, see \cite[Theorem~2]{GoMe}.
\begin{theorem}
We have
$$\liminf_{N\rightarrow \infty} \frac{E(s_n,N)}{\sqrt{N}} \geq  \frac{\sqrt{2}}{2},
$$
with probability one
with respect to the probability measure~\eqref{eq:prob_space}.
\end{theorem}
One can estimate the $N$th expansion complexity $E(s_n,N)$ in terms of the $N$th linear complexity $L(s_n,N)$, see  \cite[Theorem~3]{meniwi}.
\begin{theorem}\label{thm:MeWiNi}
Let $(s_n)$ be a sequence over $\F_q$ and let $G(x)$ bet its generating function. For $N\geq 2$, assume, that  
$$
G(x)\not\equiv 0\bmod x^N.
$$
Let $L_N=L(s_n,N)$ be the $N$th linear complexity
and let 
$$
\sum_{\ell=t_N}^{L_N} c_\ell s_{i+\ell}=0,\quad 0\le i\le N-L_N-1, 
$$
be a shortest linear recurrence for the first $N$ terms of $(s_n)$, where $c_{L_N}=1$ and $c_{t_N}\ne 0$.
 Then 
 $$
 E(s_n,N)\ge \left\{\begin{array}{ll} L_N-t_N+1 & \mbox{for } N>(L_N-t_N)(L_N-\min\{1,t_N-1\}),\\
                       \left\lceil \frac{N}{L_N-\min\{1,t_N-1\}}\right\rceil & \mbox{otherwise,}
                       \end{array}\right.
                       $$
and                       
$$
E(s_n, N)\le \min\{L_N+\max\{-1,-t_N+1\},N-L_N+2\}.
$$ 
\end{theorem}

The result formulates in a qualitative way that \emph{very} large $N$th linear complexity,  that is $N$th linear complexity close to $N$, is a non-randomness property. Moreover, it enables us to estimate the $N$th expansion complexity from below if the $N$th linear complexity is not too close to either $0$ or $N$ (in a logarithmic scale), say, of order of magnitude $\sqrt{N}$.

We refer to \cite{meniwi, howi} for applications of Theorem~\ref{thm:MeWiNi} for estimating the $N$th expansion complexity of certain sequences.

Certain subsequences of automatic sequences, say, the Thue-Morse and Rudin-Shapiro sequences along squares are not automatic, see Section~\ref{sec:normality} below, and thus have unbounded expansion complexity profile. However, their growth rates
are not known.  
For example, one can study further the Thue-Morse and Rudin-Shapiro sequence along squares.

\begin{figure}[ht]
\begin{center}
\includegraphics[scale=.51]{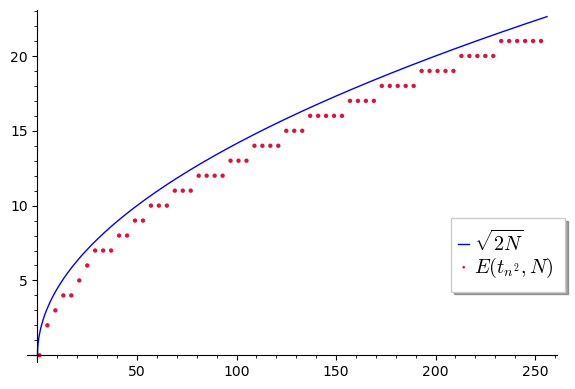}
\includegraphics[scale=.51]{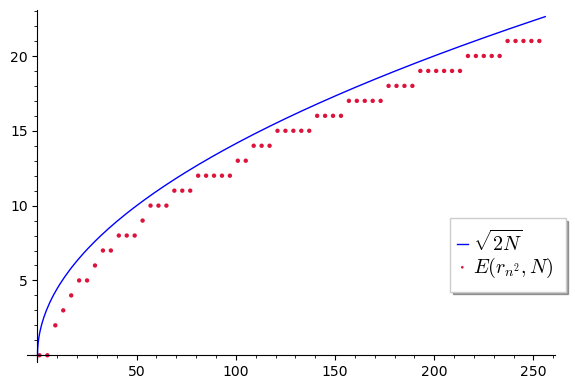} 
\end{center}
\caption{The $N$th expansion complexity of the Thue-Morse (left) and Rudin-Shapiro (right) sequence along squares.}
\label{fig:exp}
\end{figure}

\begin{problem}
Estimate the expansion complexity profiles of the subsequences $(t_{n^2})$ and $(r_{n^2})$ of the Thue-Morse and Rudin-Shapiro sequence along squares.
\end{problem}
Figure~\ref{fig:exp} suggests $E(t_{n^2},N)$ and $E(r_{n^2},N)$ are both of order of magnitude $\sqrt{N}$.

Finally, we remark that in order to use the full strength of Theorem~\ref{thm:Diem} for inferring sequence elements, one needs to require the irreducibly of the polynomial $h(x,y)$ in \eqref{eq:h}. In \cite{GoMe, GoMeNi}, the authors studied this variant of the $N$th expansion complexity and the relation between these two complexity measures.

\section{Subword complexity and normality}\label{sec:normality}

The results of Section~\ref{sec:correlation} show that many automatic sequences, including Thue-Morse and Rudin-Shapiro sequence, are balanced, that is, the frequencies of the symbols are close to the expected values. However, the frequencies of longer patterns are far from uniform. This phenomenon can be made precise by the notion of \emph{subword complexity}.

\begin{definition}
For a sequence $(s_n)$ over the alphabet $\Delta$ the \emph{subword complexity} $p(s_n,k)$ is the number of distinct subsequences of length $k$.
\end{definition}
Trivially we have $1\le p(s_n,k)\le |\Delta|^k$ and for ultimately periodic sequences we have $p(s_n,k)=O(1)$.

By \cite[Corollary 10.3.2]{alsh2} the subword complexity $p(s_n,k)$ of automatic sequences $(s_n)$ is of order of magnitude $k$.
\begin{theorem}\label{thm:subword} If $(s_n)$ is an automatic sequence that is not ultimately periodic,
then we have
\footnote{$f(k)=\Theta(g(k))$ is equivalent to $c_1g(k)\le f(k)\le c_2g(k)$ for some constants $c_2\ge c_1>0$.}
$$p(s_n,k)=\Theta(k).$$
\end{theorem}
For the Thue-Morse sequence $(t_n)$, the exact value of its subword complexity $p(t_n,k)$ was independently determined by Brlek \cite[Proposition~4.4]{br} and by de Luca and Varricchio \cite[Proposition~4.4]{deva}, see also \cite[Exercise~10.11.10]{alsh2}.
De Luca and Varricchio \cite[Property~3.3]{deva} also showed that patterns such as $000$ and $111$ do not appear in the Thue-Morse sequence
and more general the following result. 
\begin{theorem}
  The Thue-Morse sequence is cube-free, that is, no pattern of the form $www$ with $w\in \{0,1\}^k$ for some $k\ge 1$ appears in the sequence. 
\end{theorem}
The papers \cite{br,deva} contain also several other results on the non-existence of certain patterns in
the Thue-Morse sequence.

The subword complexity and the correlation measure of order $\ell$ are related by the following 
result of Cassaigne et al. \cite[Theorem 6]{cafe}.
\begin{theorem}
If for some positive integers $k$ and $N$
$$
C_{\ell}(s_n,N)\le \frac{N}{2^{2k+1}},\quad \ell=1,2,\ldots,k,
$$
then
$$p(s_n,k)=2^k.$$
\end{theorem}
For automatic sequences we can have $p(s_n,k)=2^k$ only for finitely many $k$ since $p(s_n,k)=\Theta(k)$.
However, certain subsequences of automatic squences are \emph{normal}, that is, all patterns appear in the sequence with the expected frequencies. More formally, a sequence $(s_n)$ is called \emph{normal} if for any fixed length $k$ and any pattern $\mathbf{e}\in \Delta^k$
$$
N_k(s_n, \mathbf{e}, N) =\frac{
\#\{0\le n< N:  (s_{n},s_{n+1},\ldots,s_{n+k-1})=\mathbf{e}\}}{N}\rightarrow \frac{1}{|\Delta|^k} \quad \text{as } N \rightarrow \infty.
$$

Drmota et al.\ \cite{drmari} and M\"ullner \cite{mu} proved the normality of the Thue-Morse and the Rudin-Shapiro sequences along squares, that is
\begin{equation}\label{eq:normality}
\lim_{N\rightarrow \infty}N_k(t_{n^2}, \mathbf{e}, N) =2^{-k}      
\quad \text{and} \quad
\lim_{N\rightarrow \infty}N_k(t_{n^2}, \mathbf{e}, N)=2^{-k}
\end{equation}
for any $\mathbf{e}\in \{0,1\}^k$.
The main tool to obtain the results \eqref{eq:normality} is to prove estimates on the sums $$
\sum_{n<N}(-1)^{e_0t_{n^2}+\dots + e_{k-1}t_{(n+k-1)^2} }
\quad \text{and} \quad
\sum_{n<N}(-1)^{e_0r_{n^2}+\dots + e_{k-1}r_{(n+k-1)^2} }
$$
for any $e_0,\dots, e_{k-1}\in\{0,1\}$.
These sums can be estimated via a Fourier analytic method of Mauduit and Rivat which has its origin in \cite{mari1,mari2}.
For more details we refer to the survey \cite{dr} of Drmota and the original papers \cite{drmari,mu}.

In particular, the normality results \eqref{eq:normality} yield the the subword complexities
\begin{equation}\label{eq:subword}
p(t_{n^2},k)=p(r_{n^2},k)=2^k.
\end{equation}

It is conjectured but not proved yet that the subsequences of the Thue-Morse sequence~$(t_{f(n)})$ and Rudin-Shapiro sequence~$(r_{f(n)})$ along any polynomial~$f$ of degree $d\ge 3$ 
are normal, see \cite[Conjecture 1]{drmari}.
Even the weaker problem of determining the frequency of $0$ and $1$ in the subsequences $(t_{f(n)})$ and $(r_{f(n)})$  along any polynomial $f(x)$ of degree $d\ge 3$ with $f(\N_0)\subset \N_0$ seems to be very intricate, see 
\cite[above Conjecture 1]{drmari}.

\begin{problem}
  Show that the subsequences of Thue-Morse and Rudin-Shapiro sequence along cubes, bi-squares, ..., any polynomial values for a polynomial of degree at least~$3$ are normal.  
\end{problem}

However, Moshe \cite{mo} proved the following lower bound on the subword complexity of~$(t_{f(n)})$,
\begin{equation}\label{moshe}
p(t_{f(n)},k)\ge 2^{k/2^{d-2}}.
\end{equation}
Stoll \cite{st12,st16} showed that the number of zeros (resp.\ ones) among the first $N$ sequence elements of both, $(t_{f(n)})$ and $(r_{f(n)})$,
is at least of order of magnitude $N^{4/(3d+1)}$, $d\ge 3$.
For subsequences $(z_{f(n)})$ of the Zeckendorf sum of digits sequence $(z_n)$ defined by \eqref{zeck} the numbers of zeros and ones among the first~$N$ sequence elements are both lower bounded by~$N^{4/(6d+1)}$, see Stoll~\cite{st13}.

M\"ullner and Spiegelhofer \cite{musp,sp} addressed the normality problem for the Thue-Morse sequence along the Piateski-Shapiro sequence $\lfloor n^c\rfloor$ for $1<c<3/2$. Moreover, it is asymptotically balanced (or simply normal)\cite[Theorem~1.2]{sp20} for $1<c<2$.  
For results on the Thue-Morse and Rudin-Shapiro sequence along primes see \cite{bo1,bo2,mari2,mari} and references therein.
In particular, the Thue-Morse sequence $(t_p)$ along primes is balanced, see Mauduit and Rivat \cite{mari2}. However, it is not known whether $(t_{f(p)})_{p}$
is normal for any nonconstant polynomial~$f$.

From Theorem~\ref{thm:subword} and \eqref{eq:subword} we know that $(t_{n^2})$ and $(r_{n^2})$ are not automatic
and by Theorem~\ref{thm:christol} these subsequences are, in contrast to the original sequence, not of bounded expansion complexity, that is, 
$$
\lim\limits_{N\rightarrow \infty}E(t_{n^2},N)=\lim\limits_{N\rightarrow \infty}E(r_{n^2},N)=\infty.
$$
Theorem \ref{thm:subword} combined with \eqref{moshe} implies that $(t_{f(n)})$ is not automatic
and 
$$\lim\limits_{N\rightarrow \infty}E(t_{f(n)})=\infty$$
for any polynomial of degree at least $2$ with $f(\N_0)\subset \N_0$.
Note that it was shown in \cite{al82} that $(t_{f(n)})$ is not $2^r$-automatic and in 
\cite{alsa} that $(r_{f(n)})$ is not $2^r$-automatic and thus we also have
$$
\lim\limits_{N\rightarrow \infty}E(r_{f(n)})=\infty.
$$

Subsequences of the Thue-Morse sequence along geometric sequences such as $(t_{3^n})$ seem to be even more difficult to analyze. 
For example, Lagarias \cite[Conjecture~1.12]{la} conjectured that
each pattern appears at least once in $(t_{3^n})$.
For other related results see\cite{duwe,kast}.

For more details on the normality of automatic sequences and their subsequences we refer to \cite{dr}.

\section{Analogs for finite fields}\label{sec:finite_fields}

An analog for finite fields of the problem on the distribution of automatic sequences and their subsequences was introduced by Dartyge and S\'ark\"ozy \cite{dasa}. It has been further investigated in \cite{damewi,mawi,sw18,sw1,ma}, see also
\cite{dielsh,ga17,sw2,damasa,os}.

In the finite field setting some problems can be solved although the analog for integers seems to be out of reach including the normality problem for the analog of the Thue-Morse sequence and the frequency problem for the analog of the Rudin-Shapiro sequence both along polynomials. Hence, these analogs for finite fields are further attractive sources of pseudorandomness.

For a prime $p$ and $q=p^r$ with $r\geq 2$ 
let $(\beta_1,\ldots,\beta_r)$ be an ordered basis of $\F_q$ over~$\F_p$. Then one can write all elements $\xi \in \F_q$ as
\begin{equation}\label{eq:xi}
\xi=\sum_{i=1}^r x_i\beta_i , \quad x_1,\dots, x_r\in \F_p. 
\end{equation}
It is natural to consider the coefficients $x_1,\dots, x_r$ as digits with respect to the basis~$(\beta_1,\dots, \beta_r)$. 
Then, in analogy to the Thue-Morse and Rudin-Shapiro sequence satisfying \eqref{sumofdigitsdef} we define the {\em Thue-Morse function}
$$
T\left(\sum_{i=1}^r x_i\beta_i\right)=
\sum_{i=1}^rx_i
,\quad x_1,\ldots,x_r\in \F_p,
$$
and {\em Rudin-Shapiro function}
$$
R\left(\sum_{i=1}^r x_i\beta_i\right)
=\sum_{i=1}^{r-1}x_ix_{i+1}
,\quad x_1,\ldots,x_r\in \F_p,
$$
on $\F_q$.

Dartyge and S\'ark\"ozy \cite{dasa} studied the balance of the Thue-Morse function along polynomial values. They derived results using the Weil bound \cite[Theorem~5.38]{lini} on additive character sums:
\begin{lemma}\label{lemma:weil}
Let $f\in\F_q[x]$ be of degree $d\geq 1$ with $\gcd(d,q)=1$ and $\psi$ be a nontrivial additive character of $\F_q$. Then
$$
\left|\sum_{\xi\in\F_q}\psi(f(\xi))\right|\leq (d-1)\sqrt{q}.
$$
\end{lemma}
Put
$$
e(\alpha)=\exp(2\pi i \alpha), \quad \alpha\in\mathbb{R},
$$
and 
note that $\psi(x)=e(T(x)/p)$ is a nontrivial additive character of $\F_q$.
Then from
$$
\sum_{h=0}^{p-1} e\left(\frac{ha}{p}\right)=\left\{\begin{array}{cc}0,& a\neq 0,\\ p, &a=0,\end{array}\right.
\quad a\in \F_p,
$$
we get
$$
\#\{\xi\in \F_q:  T(f(\xi))=c\}= \frac{1}{p}\sum_{h=0}^{p-1}\sum_{\xi\in \F_q}\psi\left(hf(\xi)\right)e\left(\frac{-hc}{p}\right).
$$
The contribution of $h=0$ is trivially $p^{r-1}$, which is the expected  number of solutions.  The other terms for $h\neq 0$ contribute to the error term and can be bounded by Lemma~\ref{lemma:weil}. We immediately get \cite[Theorem~1.2]{dasa}:

\begin{theorem}\label{thm:dasa}
Let $f\in \F_q[x]$ be of degree $d$ with $\gcd(d, q) = 1$. Then for all $c\in \F_p$, we
have
$$
\big|\#\{\xi\in \F_q:  T(f(\xi))=c\}-p^{r-1}\big|\le (d-1)p^{r/2}.
$$
\end{theorem}

Later Dartyge, M\'erai and Winterhof \cite{damewi} investigated this problem for the Rudin-Shapiro function. The main difference between the two problems is that the Rudin-Shapiro function is not a linear map contrary to the Thue-Morse function. Standard character sum techniques fail in this situation. Namely, consider $R(f(\xi))$ with $\xi$ having the form \eqref{eq:xi} as a polynomial in the $r$ variables $x_1,\dots, x_r$. Then using Lemma~\ref{lemma:weil} for one coordinate $x_i$ one gets an error term larger than the main term. Stronger results in higher dimension such as the Deligne bound \cite[Th\'eor\`eme 8.4]{De74} also cannot be applied as it needs some more  technically intricate conditions which are not satisfied in our situation. However, sacrificing the explicit dependence of the degree $d$, one can use an affine version of the Hooley-Katz Theorem, see \cite{ho} or \cite[Theorem~7.1.14]{handbookFF}.

First recall that the \emph{(affine) singular locus} $\cL(F)$ of a polynomial $F\in\F_p[x_1,\dots, x_r]$ is the set of common zeros in $\overline{\F_p}^r$
of the polynomials\footnote{$\overline{\F_p}=\bigcup_{n=1}^\infty \F_{p^n}$ denotes the algebraic closure of $\F_p$.}
$$
F,\frac{\partial F}{\partial x_1},\ldots,\frac{\partial F}{\partial x_r}.
$$
We also recall that the dimension of $\cL(F)$ is the largest $d$ for which there exist $1\le i_1<i_2<\ldots<i_d\le r$ such that there is no nonzero polynomial $P$ in $d$ variables with
$P(y_{i_1},\ldots,y_{i_d})=0$ for all $(y_1,\ldots,y_r)\in \cL(F)$, see \cite[Corollary~9.5.4]{ideal}.

\begin{lemma}\label{lemma:HK}
Let $Q\in\F_p[x_1,\dots, x_r]$ be of degree~$d\ge 1$ such that the dimensions of the singular loci of $Q$ and its homogeneous part $Q_d$ of degree $d$ satisfy
$$
\max\{\dim(\cL(Q)),\dim(\cL(Q_d))-1\}\le s.
$$ 
Then  the number $N$ of zeros of~$Q$ in $\F_p^r$ satisfies
$$
\left|N-p^{r-1}\right|\le C_{d,r}p^{(r+s)/2},
$$
where $C_{d,r}$ is a constant depending only on $d$ and $r$.
\end{lemma}
Then using Lemma~\ref{lemma:HK}, one can show that the Rudin-Shapiro function is also asymptotically balanced on polynomial values, see \cite[Theorem~1]{damewi}.
\begin{theorem}\label{thm:damewi}
Let $f\in \F_q[x]$ be of degree $d$ with $\gcd(d, q) = 1$. Then for all $c\in \F_p$, we
have
$$
\big|\#\{\xi\in \F_q:  R(f(\xi))=c\}-p^{r-1}\big|\le C_{d,r} p^{(3r+1)/4},
$$
where the constant $C_{d,r}$ depends only on the degree $d$ of $f$ and $r$.
\end{theorem}

 Theorem \ref{thm:damewi} is nontrivial if $r$ is fixed and $p\rightarrow \infty$. Contrary to Theorem \ref{thm:dasa}, nothing is known for the dual situation.
\begin{problem}\label{prob:RS_FF}
For fixed prime $p$ show that if $r$ is large enough, then the Rudin-Shapiro function along polynomial values is balanced possibly under some natural restrictions on the polynomial.
\end{problem}

Analogously to the normality results of Section~\ref{sec:normality}, Makhul and Winterhof \cite{mawi} obtained results on the normality of the Thue-Morse function along polynomial values. For sake of simplicity we state the case when the polynomial $f$ has degree $d$ smaller than the characteristic $p$, \cite[Corollary~1]{mawi}.

\begin{theorem}
Assume $1\leq  d < p$ and $s\leq  d$. For any polynomial $f\in\F_q[x]$ of degree $d$ and any pairwise distinct $\alpha_1, \dots, \alpha_s\in\F_q$  and any $c_1,\dots, c_s\in \F_p$ we have
$$
\big|\#\{\xi\in \F_q:  T(f(\xi+\alpha_i))=c_i, 1\leq i\leq s\}-p^{r-s}\big|\le (d-1) p^{r/2}.
$$
\end{theorem}
Note that the restriction $s\le d$ is natural and counterexamples for $s>d$ are easy to construct. 

The case of the Rudin-Shapiro function is much more intricate.

\begin{problem}
Study the normality of the Rudin-Shapiro function at $f(x)$. Namely, show that
$$
\frac{\#\{\xi\in \F_q:  R(f(\xi+\alpha_i))=c_i, 1\leq i\leq s\}}{p^{r-s}} \rightarrow 1 \quad \text{as } p\rightarrow \infty
$$
for some $s\geq 2$ and any $f\in\F_q[x]$ of fixed degree.
\end{problem}

Of course, this problem is also open for fixed $p$ and $r\rightarrow \infty$ even in the simplest 
case~$s=1$, see Problem~\ref{prob:RS_FF}.

It is natural to define the 
{\em Rudin-Shapiro function} on the polynomial ring~$\F_p[t]$ by assigning the coefficients of the polynomial $f(t)\in\F_p[t]$ to $(x_1,\dots, x_r)$, that is,
$$
R(t^r+x_{1}t^{r-1}+\dots + x_r )=\sum_{i=1}^{r-1}x_ix_{i+1},
$$
for $x_1,\dots, x_r\in \F_p$.

Analogously to the result of Mauduit and Rivat \cite{mari} on the Rudin-Shapiro sequence along prime numbers, it is natural to investigate the balance and the normality of the Rudin-Shapiro function along irreducible polynomials. As the number of monic irreducible polynomials of degree $r$ is $p^r/r+o(p^r)$, see for example \cite[Theorem~3.25]{lini},
we expect that the frequency of each element $c$ is $\frac{p^{r-1}}{r}+o(p^{r-1})$.
For $r=2$ and fixed $c\in \F_p$ we have to count the number of $x_2\in \F_p^*$ such that $t^2+x_2^{-1}ct+x_2$ is irreducible over $\F_p$ or equivalently the discriminant $x_2^{-2}c^2-4x_2$ is a quadratic non-residue modulo~$p$. This number is
$$\frac{1}{2}\sum_{x_2\in \F_p^*} \left(1-\left(\frac{c^2-4x_2^3}{p}\right)\right)=
\left\{\begin{array}{cc} \frac{p-1}{2}, & c=0,\\
\frac{p-1}{2}+O(p^{1/2}), & c\ne 0,\end{array}\right.$$
by the Weil bound for multiplicative character sums \cite[Theorem 5.41]{lini}, where~$\left(\frac{.}{.}\right)$ is the Legendre symbol.
\begin{problem}
Prove that for all $c\in\F_p$ and $r\ge 3$ we have
$$
\lim_{p\rightarrow\infty}\frac{ \#\{f\in \F_p[t]: \deg f=r, f \text{ monic and irreducible over $\F_p$}, R(f)=c \} }{p^{r-1}}=\frac{1}{r}.$$
\end{problem}

We remark that one can define the {\em Thue-Morse function}
by
$$T(f)=T(t^r+x_1t^{r-1}+\ldots+x_r)=x_1+\ldots+x_r=f(1)-1.$$
Note that for irreducible polynomials $f(x)$ we have $T(f)\ne -1$ and for $c\ne -1$ the number of monic irreducible polynomials of degree $r=2$
with $T(f)=c$ is
$$
\frac{1}{2}\sum_{u\in \F_p\atop u^2\ne c+1}\left(1-\left(\frac{u^2-c-1}{p}\right)\right)  =\frac{p-\left(\frac{c+1}{p}\right)}{2},
$$
where we used a well-known result on sums of Legendre symbols of quadratic polynomials, see for example \cite[Theorem~5.48]{lini}.
In general, since $f(x)$ is irreducible whenever $f(x-1)$ is irreducible we have to estimate the number $I_c$ of monic irreducible polynomials with fixed constant term $c\ne 0$ which satisfies
$$\frac{1}{r}\left(\frac{p^r-1}{p-1}-2p^{r/2}\right)\le I_c\le \frac{p^r-1}{r(p-1)},$$
see \cite{car} or \cite[Theorem~3.5.9]{handbookFF}, and we get the desired 
$$
I_c=\frac{p^{r-1}}{r}+o(p^{r-1})
$$
for $r\ge 3$ as well.

Moreover, the corresponding normality problem is trivial since for any polynomial~$g(x)$ of degree at most $r-1$ the value $T(f+g)=f(1)+g(1)-1$ is uniquely defined by $T(f)=f(1)-1$ and $g(1)$.

For other results on 'digits' along irreducible polynomials see for example
\cite[Chapter~3]{handbookFF} and
\cite{gakuwa,gr,pol,tuwa,ha,por}.

\section*{Acknowledgment}
The authors were supported by the Austrian Science Fund FWF grants P 30405 and P~31762.

They wish to thank Jean-Paul Allouche, Harald Niederreiter, Igor Shparlinski, Cathy Swaenepoel, Thomas Stoll and Steven Wang for very useful discussions.

\end{document}